%% file: ex_article.tex
\documentclass[review,hidelinks,onefignum,onetabnum]{siamart251216}
\usepackage{silence}
\input{ex_shared}

\usepackage{amsmath}
\usepackage{amssymb}
\usepackage{microtype}
\usepackage{subcaption}
\usepackage{graphicx}
\usepackage{booktabs}
\usepackage{multirow}
\usepackage{tabularx} 
\usepackage{verbatim}
\usepackage{hyperref}
\usepackage{cleveref}
\usepackage{bm}
\usepackage{enumitem}
\usepackage{xcolor}

\newif\ifshowrevisions
\showrevisionstrue
\newif\ifshowauthorchecks
\showauthorcheckstrue

\newcommand{\eps}{\epsilon}

\newcommand{\bA}{\mathbf{A}}

\newcommand{\bQ}{\mathbf{Q}}
\newcommand{\bLambda}{\mathbf{\Lambda}}
\newcommand{\curl}{\operatorname{curl}}

\newcommand{\dx}{\, dx}

\newcommand{\inner}[2]{\left\langle#1, #2\right\rangle}

\crefname{assumption}{assumption}{assumptions}
\Crefname{assumption}{Assumption}{Assumptions}

\ifpdf
\hypersetup{
  pdftitle={Spectral Gaps and Low-Energy Modes in Anisotropic High-Contrast PDEs},
  pdfauthor={Eric T. Chung, Lijian Jiang and Yajun Wang}
}
\fi

\begin{document}

\maketitle

\begin{abstract}
We study the low-energy spectrum of high-contrast partial differential equations with strongly anisotropic inclusions.
The inclusion stiffness and mass are scaled by a contrast parameter \(\eta\), while a small material eigenvalue \(\epsilon\) controls anisotropy.
Stable extensions of local kernel fields give eigenvalue upper bounds proportional to the inverse contrast, uniformly in anisotropy.
Additional low-energy modes arise from fields whose nonzero differential images lie in weak material directions.
Each fixed family of compatible weak fields yields additional eigenvalue upper bounds of order $\epsilon+ \eta^{-1}$, provided a stable extension exists.
Under discrete stability assumptions, we obtain corresponding finite element bounds for the weak modes represented by the mesh.
For separable scalar diffusion, we prove norm-resolvent convergence to independent background Dirichlet and inclusion Neumann operators as contrast tends to infinity with anisotropy fixed.
Numerical experiments for four PDE models illustrate weak modes and mesh-dependent spectral transitions, and scalar diffusion experiments also test spectral decoupling and eigenvalue counting.
\end{abstract}



\section{Introduction}
High-contrast partial differential equations arise in composite materials, porous media, and electromagnetic structures. Their spectra reflect both material properties and the geometry of heterogeneous phases. Classical homogenization and two-scale convergence provide fundamental tools for studying effective behavior \cite{CioranescuDonato1999,Allaire1992}, including elliptic eigenvalue problems \cite{Kesavan1979}. Homogenization and resolvent analyses provide tools for studying the resulting limiting operators and their spectral behavior \cite{Zhikov2000,KamotskiSmyshlyaev2019,CherednichenkoCooper2016}. In spectral coarse-space constructions, the number of low-energy modes guides the choice of coarse-space dimension.

Related computational developments include multiscale finite element methods \cite{HouWu1997}, variational multiscale formulations \cite{HughesEtAl1998}, heterogeneous multiscale methods \cite{AbdulleEtAl2012}, and subsequent multiscale finite element theory \cite{EfendievHou2009}. Localized orthogonal decomposition further extends operator-adapted approximation \cite{MalqvistPeterseim2014,AltmannHenningPeterseim2021}. Spectral information enters more explicitly in generalized multiscale finite element methods \cite{EfendievGalvisHou2013,ChungEfendievFu2014,chung2023multiscale}, constraint energy minimizing constructions \cite{ChungEfendievLeung2018,chung2023multiscale}, optimal local approximation spaces \cite{BabuskaLipton2011,MaScheichlDodwell2022}, and localized subspace iteration \cite{GuanEtAl2026}. Spectral selection also supports robust domain decomposition through high-contrast coarse spaces \cite{GalvisEfendiev2010,EfendievEtAl2012}, generalized eigenproblems in overlaps \cite{SpillaneEtAl2014}, and adaptive FETI-DP and BDDC constraints \cite{KlawonnEtAl2015,KimChungWang2017}. These methods motivate the question studied here: how does the number of low-energy modes depend on contrast, anisotropy, and mesh resolution?

Local operator kernels provide a basic source of small eigenvalues in high-contrast media. Examples include constants for diffusion, rigid motions for elasticity, and affine functions for biharmonic problems. For problems associated with differential complexes, such as Maxwell's equations, the analysis must also account for the kernel structure and compatibility conditions \cite{ArnoldFalkWinther2006,Monk2003}.
Our preceding work \cite{ChungJiangWangTopological} related local null spaces and topological structure to spectral separation. 
Under suitable coercivity assumptions on the orthogonal complement of the kernel, the kernel-generated eigenvalues are bounded by $C\eta^{-1}$, while the rest of the spectrum stays uniformly bounded away from zero as the contrast $\eta\to\infty$.

The present paper investigates how strong anisotropy affects the low-energy spectrum and the spectral separation between kernel-generated modes and additional modes. We consider separated inclusions whose stiffness and mass are scaled by $\eta\gg1$, with a weak material eigenvalue controlled by $\epsilon\ll1$.  
We call trial fields whose differential images lie in the weakly penalized material subspace compatible weak fields. Such fields can have small Rayleigh quotients even when they do not belong to the local kernel, giving rise to additional low-energy modes. Moreover, counting eigenvalues in a fixed energy window as $\epsilon\to0$ differs from classical high-energy asymptotics \cite{SafarovVassiliev1997}.

Our analysis connects three aspects of the low-energy spectrum. Stable kernel extensions and compatible weak fields give the bounds \(\lambda^{\eta,\epsilon}_N\le C\eta^{-1}\) and \(\lambda^{\eta,\epsilon}_{N+K}\le C_K(\epsilon+\eta^{-1})\) for each fixed \(K\). Under suitable assumptions on discrete liftings and reduced spectra, we derive finite element eigenvalue bounds for resolved weak modes. Finally, for separable scalar diffusion, we prove norm-resolvent convergence at fixed \(\epsilon>0\) to a direct sum of background Dirichlet and inclusion Neumann operators, and derive an exact counting formula for this limit. The numerical examples relate observed spectral transitions to the combined dimensions of local discrete kernel and weak spaces. Additional scalar experiments examine spectral decoupling and compare global finite element eigenvalue counts with exact counts and asymptotics for the decoupled limit.

The remainder of the paper is organized as follows. \Cref{sec:pro_setting} introduces the high-contrast anisotropic eigenvalue problem and its variational formulation. \Cref{sec:kernel-modes} establishes upper bounds for kernel-generated eigenvalues, while \Cref{sec:aniso_spectrum} derives additional low-energy bounds from compatible weak fields. \Cref{sec:instantiations} provides concrete constructions for diffusion, elasticity, biharmonic, and Maxwell models. \Cref{sec:fem-truncation} studies finite element truncation under discrete lifting and reduced-spectrum assumptions. \Cref{sec:modified-weyl} establishes norm-resolvent decoupling for scalar diffusion at fixed $\epsilon>0$ as $\eta\to\infty$ and derives exact counting formulas and anisotropic asymptotics for the decoupled limit. \Cref{sec:numex} presents numerical experiments on weak modes, mesh-dependent spectral transitions, scalar spectral decoupling, and eigenvalue counting. \Cref{sec:conclu} summarizes the main findings and discusses future directions.

\section{High-contrast anisotropic spectral problem}
\label{sec:pro_setting}
Let \(\Omega\subset\mathbb R^d\), \(d=2,3\), be a bounded Lipschitz domain. Let \(D_1,\ldots,D_M\) be connected bounded Lipschitz inclusions whose closures are mutually disjoint and compactly contained in \(\Omega\). The geometry is fixed independently of \(\eta\) and \(\epsilon\). 
Set \(D=\bigcup_{i=1}^M D_i\) and assume that the background \(\Omega_0=\Omega\setminus\overline D\) is connected. Additional topological assumptions for the Maxwell problem are specified in Section~\ref{sec:maxwell-fields}.

Let $V(\Omega)$ denote an appropriate global Hilbert space defined over $\Omega$ enforcing homogeneous boundary conditions on $\partial \Omega$. Depending on the physical context, $V(\Omega)$ will be instantiated as $H_0^1(\Omega)$, $[H_0^1(\Omega)]^d$, $H_0^2(\Omega)$, or $H_0(\curl, \Omega)$. 
We write \(\mathcal L:V(\Omega)\to L^2(\Omega;\mathbb E)\), where \(\mathbb E\) is the finite-dimensional Euclidean space of differential images, with \(s=\dim\mathbb E\). Let \(m\) denote the differential order of \(\mathcal L\). For symmetric tensor fields, the inner product on \(\mathbb E\) is the Frobenius inner product. The local spaces \(V(D_i)\) are defined without essential boundary conditions on \(\partial D_i\). They are \(H^1(D_i)\), \([H^1(D_i)]^d\), and \(H^2(D_i)\) for the diffusion, elasticity, and biharmonic models, respectively. The constrained spaces used for Maxwell fields are defined separately in Section~\ref{sec:maxwell-fields}.

We investigate the spectral properties of the generalized eigenvalue problem: Find eigenpairs $(\lambda^{\eta,\eps}, u^{\eta,\eps}) \in (0,\infty) \times V(\Omega)$, with $u^{\eta,\eps}\ne0$, such that
\begin{equation}\label{eq:abstract_eigen}
    a^{\eta,\eps}(u^{\eta,\eps}, v) = \lambda^{\eta,\eps} m^\eta(u^{\eta,\eps}, v), \quad \forall v \in V(\Omega).
\end{equation}

The multiscale nature of the problem is encoded in the parameter $\eta \gg 1$, which represents the extreme contrast ratio between the inclusions and the background matrix. The stiffness and mass bilinear forms are defined as follows:
\begin{equation}
\begin{aligned}
    a^{\eta,\eps}(u, v) &= \int_{\Omega_0} \inner{\mathcal{L} u}{\mathcal{L} v} \dx + \eta \sum_{i=1}^M \int_{D_i} \inner{\bA^{\eps}(x) \mathcal{L} u}{\mathcal{L} v} \dx, \\
    m^\eta(u, v) &= \int_{\Omega_0} u \cdot v \dx + \eta \sum_{i=1}^M \int_{D_i} u \cdot v \dx.
\end{aligned}
\end{equation}

The mass density is \(\rho_\eta=\eta\) in \(D\) and \(\rho_\eta=1\) in \(\Omega_0\), so \(m^\eta(u,v)=\int_\Omega\rho_\eta u\cdot v\,dx\). On a fixed finite-dimensional trial space, a local kernel field contributes no inclusion energy, while its inclusion mass is proportional to \(\eta\). A background extension whose energy is bounded independently of \(\eta\) therefore yields an \(O(\eta^{-1})\) upper bound for the Rayleigh quotient. For a compatible weak field, the inclusion energy scales with \(\eta\epsilon\), yielding an \(O(\epsilon)\) contribution to the Rayleigh quotient. A background extension with uniformly bounded energy contributes an additional \(O(\eta^{-1})\) term.

For almost every \(x\in D\), the material tensor \(\bA^\eps(x):\mathbb E\to\mathbb E\) is symmetric positive definite. 
By the spectral theorem for symmetric tensors, for almost every $x \in D_i$, there exists an orthogonal matrix $\bQ(x)$ and a diagonal matrix $\bLambda(x)$ such that
\begin{equation}
    \bA^{\eps}(x) = \bQ(x) \bLambda(x) \bQ^T(x).
\end{equation}
The columns of $\bQ(x)$, denoted by $\{q_1, \dots, q_s\}$, form the orthonormal principal axes of the material, and $\bLambda(x) = \operatorname{diag}(\gamma_1, \dots, \gamma_s)$ has the corresponding eigenvalues as its diagonal entries.

To model strong anisotropy, we impose the following separation of scales between the weak and strong eigenvalues of $\bA^\eps$.

\begin{assumption}[Extreme anisotropy]
\label[assumption]{assump:extreme_anisotropy}
We assume that \(\bA^\eps(x)\) has one weak eigenvalue \(\epsilon\), while its remaining \(s-1\) eigenvalues are bounded above and below by positive constants independent of \(\eta\) and \(\epsilon\).
\begin{equation}
    C_1 \le \gamma_1, \dots, \gamma_{s-1} \le C_2, \quad \text{and} \quad \gamma_s = \eps \ll 1.
\end{equation}
\end{assumption}

Let \(P_s(x)\) denote the orthogonal projection onto \(\operatorname{span}\{q_1(x),\ldots,q_{s-1}(x)\}\). 
We call a local trial field \(w\) a compatible weak field if
\[
P_{\mathrm s}(x)\mathcal Lw(x)=0
\qquad\text{for almost every }x\in D.
\]
Equivalently, its differential image lies in the weak material eigenspace almost everywhere in \(D\).
Using this orthogonal decomposition, the local stiffness energy within an inclusion $D_i$ can be explicitly decoupled into a sum of directional energies along the principal axes:
\begin{equation}\label{eq:energy_decoupled}
    a_{D_i}(u, u) = \eta \int_{D_i} \left( \sum_{k=1}^{s-1} \gamma_k |\inner{\mathcal{L} u}{q_k}|^2 + \eps |\inner{\mathcal{L} u}{q_s}|^2 \right) \dx.
\end{equation}

Assume that the energy form \(a_{\eta,\epsilon}\) is closed on \(V(\Omega) \) and that its form domain is compactly embedded in the weighted \(L^2\) space. The eigenvalues then admit the standard min–max characterization \cite{Kato1995}. 
\begin{equation}
    R^{\eta,\eps}(u) = \frac{a^{\eta,\eps}(u, u)}{m^\eta(u, u)}, \quad \lambda_n^{\eta,\eps} = \min_{\substack{S_n \subset V(\Omega) \\ \dim(S_n)=n}} \max_{u \in S_n \setminus \{0\}} R^{\eta,\eps}(u).
\end{equation}
For Maxwell problems, the min–max characterization is understood in the weighted gradient-orthogonal space $\mathcal W_{\eta}$ introduced in Section \ref{sec:maxwell-fields}. The full space $H_0(\curl, \Omega)$ contains an infinite-dimensional gradient kernel and is not used for numbering the positive eigenvalues.

\section{Kernel-generated low-energy modes}
\label{sec:kernel-modes}
For each inclusion, let
\[
\mathcal K_i=\{v\in V(D_i):\mathcal Lv=0\},\qquad
\mathcal K=\bigoplus_{i=1}^M\mathcal K_i,\qquad
\|k\|_D^2=\sum_{i=1}^M\|k_i\|_{L^2(D_i)}^2.
\]
\begin{assumption}[Finite local kernels and stable extension]
\label[assumption]{assump_ker_and_ext}
The space \(\mathcal K\) has finite dimension \(N\). There is a linear operator \(E_0:\mathcal K\to V(\Omega)\), independent of the material parameters, such that
\[
(E_0k)|_{D_i}=k_i,\qquad
\|\mathcal LE_0k\|_{L^2(\Omega_0)}^2\le C_0\|k\|_D^2
\quad\text{for all }k\in\mathcal K.
\]
    
\end{assumption}
\begin{lemma}[Extension of polynomial kernel fields]
    \Cref{assump_ker_and_ext} holds for the local kernels of the gradient, symmetric gradient, and Hessian operators on the separated inclusions considered here.
\end{lemma}
\begin{proof}
    Their local kernels consist, respectively, of constants, infinitesimal rigid motions \cite{Ciarlet1988}, and affine functions. Each kernel field therefore has a polynomial representative on \(\mathbb R^d\). Choose smooth cutoff functions \(\chi_i\) supported in mutually disjoint neighborhoods compactly contained in \(\Omega\), with \(\chi_i\equiv1\) on a neighborhood of \(\overline{D_i}\), and set \(E_0k=\sum_i\chi_i k_i\). The extension \(E_0k\) agrees with \(k_i\) on each inclusion \(D_i\) and satisfies the homogeneous boundary conditions on \(\partial\Omega\). 
    Since \(\mathcal K\) is finite dimensional, the fixed linear map \(k\mapsto(\mathcal L E_0k)|_{\Omega_0}\) is bounded from \((\mathcal K,\|\cdot\|_D)\) into \(L^2(\Omega_0;\mathbb E)\). This yields the stated background estimate with a constant independent of \(\eta\) and \(\epsilon\).
\end{proof}

\begin{theorem}[Uniform upper bound for kernel-generated modes]
Under \Cref{assump_ker_and_ext}, if \(N\ge1\),
\[
\lambda^{\eta,\epsilon}_N\le C_0\eta^{-1}.
\]
\end{theorem}
\begin{proof}
Since \((E_0k)|_{D_i}=k_i\) for every \(i\), the extension operator \(E_0\) is injective, and hence \(\dim(E_0\mathcal K)=\dim\mathcal K=N\). For every \(k\in\mathcal K\), the inclusion energy vanishes. Therefore
\[
a^{\eta,\epsilon}(E_0k,E_0k)\le C_0\|k\|_D^2, \qquad m^\eta(E_0k,E_0k)\ge\eta\|k\|_D^2.
\]
Applying the min--max principle to the \(N\)-dimensional trial space \(E_0\mathcal K\) yields the stated bound.
\end{proof}

This estimate guarantees at least \(N\) eigenvalues with an \(O(\eta^{-1})\) upper bound, counted with multiplicity. It does not exclude additional eigenvalues on the same scale.
The integer \(N\) denotes the total dimension of the local kernels. Its interpretation depends on the differential operator.
Maxwell fields require the constrained construction in Section~\ref{sec:maxwell-fields}.

\section{Anisotropy-induced weak modes}
\label{sec:aniso_spectrum}
We enlarge the local kernel by compatible weak fields and extend the resulting trial space to the global domain.

\begin{theorem}
\label{thm:anisotropy_spectrum}
Suppose \Cref{assump:extreme_anisotropy,assump_ker_and_ext} hold. For every fixed integer \(K\ge1\), assume there is a \(K\)-dimensional space
\[
\mathcal W_K\subset\left(\bigoplus_{i=1}^M V(D_i)\right)\cap\mathcal K^{\perp_{L^2(D)}}
\]
such that
\[
P_{\mathrm s}\mathcal Lw=0\quad\text{in }D,
\qquad
\|\mathcal Lw\|_{L^2(D)}^2\le C_K\|w\|_D^2,
\qquad  \forall w\in\mathcal W_K.
\]
Assume also that there is a linear extension \(E_K:\mathcal K\oplus\mathcal W_K\to V(\Omega)\) satisfying
\[
(E_Kz)|_{D_i}=z_i,\qquad
\|\mathcal LE_Kz\|_{L^2(\Omega_0)}^2\le C_K\|z\|_D^2,
\qquad z\in\mathcal K\oplus\mathcal W_K.
\]
The spaces and constants may depend on \(K\) and the fixed geometry, but the constants are independent of \(\eta\) and \(\epsilon\). Then
\[
\lambda^{\eta,\epsilon}_{N+K}\le C_K(\epsilon+\eta^{-1}).
\]
\end{theorem}
\begin{proof}
    Write \(z=k+w\), with \(k\in\mathcal K\) and \(w\in\mathcal W_K\). Local orthogonality yields \(\|z\|_D^2=\|k\|_D^2+\|w\|_D^2\). Since \(\mathcal Lk=0\) and \(\mathcal Lw\) has no strong component,
\[
\begin{aligned}
a^{\eta,\epsilon}(E_Kz,E_Kz)
&\le C_K\|z\|_D^2+\eta\epsilon\|\mathcal Lw\|_{L^2(D)}^2
\le C_K(1+\eta\epsilon)\|z\|_D^2,\\
m^\eta(E_Kz,E_Kz)&\ge\eta\|z\|_D^2.
\end{aligned}
\]
The extension is injective, so its image has dimension \(N+K\). Applying min–max on this entire trial space proves the bound.
\end{proof}

Under the stated assumptions, any fixed number of compatible weak fields yields additional small eigenvalues as \(\epsilon+\eta^{-1}\to0\).
The dependence of $C_K$ on $K$ limits this conclusion to fixed-dimensional trial spaces.
\Cref{sec:instantiations} presents concrete PDE models that satisfy the stated assumptions.

\section{Instantiations in physical PDEs}
\label{sec:instantiations}
This section verifies the assumptions of \Cref{thm:anisotropy_spectrum} for concrete diffusion, elasticity, and biharmonic models with constant material directions. We treat Maxwell fields separately using a constrained construction.

\subsection{Diffusion, elasticity, and biharmonic models}
Let \(q\) be a fixed unit vector, and consider trial fields that depend on \(x\) only through \(q\cdot x\). For elasticity, also choose a unit vector \(p\perp q\) to specify the displacement direction. \Cref{tab:sobolev-models} lists the trial fields, their differential images, and the corresponding weak material states. These states are vectors for diffusion and symmetric tensors for elasticity and the biharmonic model. In each case, the material tensor \(A^\epsilon\) has eigenvalue \(\epsilon\) along the indicated state and is uniformly positive definite on its orthogonal complement.
Here
\[
e(u)=\tfrac12(\nabla u+\nabla u^T),\qquad
S=\frac{p\otimes q+q\otimes p}{\sqrt2},\qquad R_i=\dim\mathcal K_i.
\]
The differential order is \(m=1\) for diffusion and elasticity and \(m=2\) for the biharmonic model.

\begin{table}[htbp]
\centering
\caption{Local kernel dimensions and compatible weak fields for the Sobolev models.}
\label{tab:sobolev-models}
\renewcommand{\arraystretch}{1.25}
\begin{tabularx}{\textwidth}{@{}Xlll@{}}
\toprule
\textbf{Model} & \textbf{Weak state} & \textbf{Trial field} $w$ & \textbf{Image} $\mathcal Lw$ \\
\midrule
Diffusion, $\mathcal L=\nabla$, $R_i=1$ & $q$ & $\phi(q\cdot x)$ & $\phi'(q\cdot x)q$ \\
Elasticity, $\mathcal L=e$, $R_i=d(d+1)/2$ & $S$ & $\phi(q\cdot x)p$ & $\phi'(q\cdot x)S/\sqrt2$ \\
Biharmonic, $\mathcal L=\nabla^2$, $R_i=d+1$ & $q\otimes q$ & $\phi(q\cdot x)$ & $\phi''(q\cdot x)q\otimes q$ \\
\bottomrule
\end{tabularx}
\end{table}

In every case, \(P_{\mathrm s}\mathcal Lw=0\). Smooth functions \(\phi\) of one variable generate infinitely many linearly independent local weak fields through the constructions above. Let \(\Pi_{\mathcal K}\) denote the \(L^2(D)\)-orthogonal projection onto \(\mathcal K\). Each field \(w-\Pi_{\mathcal K}w\) is orthogonal to \(\mathcal K\) and has the same differential image as \(w\), since \(\mathcal L(\Pi_{\mathcal K}w)=0\). Because \(\mathcal K\) is finite dimensional, the resulting fields still span an infinite-dimensional space. 
For each fixed \(K\), choose \(K\) linearly independent projected fields and let \(\mathcal W_K\) be their span. Bounded Sobolev extension operators \cite{AdamsFournier2003}, combined with fixed smooth cutoffs around the separated inclusions, yield a linear extension \(E_K:\mathcal K\oplus\mathcal W_K\to V(\Omega)\) that preserves the fields on each inclusion. Norm equivalence on this fixed finite-dimensional space gives the estimates required in \Cref{thm:anisotropy_spectrum}, with constants independent of \(\eta\) and \(\epsilon\).
Consequently, for each fixed \(K\), \Cref{thm:anisotropy_spectrum} yields
\[
\lambda^{\eta,\epsilon}_{N+K}\le C_K(\epsilon+\eta^{-1}).
\]
The elasticity and biharmonic constructions apply in $d=2,3$, and \(S\) represents a weak pure-shear strain state.

\subsection{Maxwell model}
\label{sec:maxwell-fields}
For the Maxwell model, set \(\mathcal L=\operatorname{curl}\), so that \(m=1\). Assume that \(\Omega\subset\mathbb R^3\) is a bounded simply connected Lipschitz domain with connected boundary. These assumptions hold for the cubic domain used in the numerical experiments.
With the density \(\rho_\eta\) defined in \Cref{sec:pro_setting}, set
\[
\mathcal W_\eta=\{u\in H_0(\operatorname{curl};\Omega):
m^\eta(u,\nabla p)=0\ \text{for all }p\in H_0^1(\Omega)\}.
\]
Every eigenfield associated with a positive eigenvalue belongs to \(\mathcal W_\eta\), as follows by testing the eigenvalue equation with \(\nabla p\). Under the stated topological assumptions, the curl kernel in \(H_0(\operatorname{curl};\Omega)\) is exactly \(\nabla H_0^1(\Omega)\), so restricting the eigenproblem to its \(m^\eta\)-orthogonal complement \(\mathcal W_\eta\) removes this kernel.

For each fixed \(1\le\eta<\infty\), the defining orthogonality of \(\mathcal W_\eta\) implies \(\operatorname{div}(\rho_\eta u)=0\) in the distributional sense. The weighted Maxwell compactness theorem \cite{BauerPaulySchomburg2016} then shows that \(\mathcal W_\eta\), equipped with the \(H(\operatorname{curl};\Omega)\) norm, is compactly embedded in the weighted \(L^2\) space. Together with the removal of the curl kernel, this yields a discrete positive spectrum and allows us to apply the min--max characterization for fixed material parameters.

On each inclusion, we introduce the local space
\[
\mathcal Z_i=\{v\in H(\operatorname{curl};D_i):
(v,\nabla q)_{D_i}=0\ \text{for all }q\in H^1(D_i)\},
\qquad
\mathcal H_i=\mathcal Z_i\cap\ker(\operatorname{curl}).
\]
The space \(\mathcal Z_i\) consists of fields that are \(L^2(D_i)\)-orthogonal to all gradients of functions in \(H^1(D_i)\). Equivalently, these fields have zero divergence in \(D_i\) and zero normal trace on \(\partial D_i\), both in the weak sense.
The curl-free subspace \(\mathcal H_i\) of \(\mathcal Z_i\) is the space of harmonic Neumann fields on \(D_i\). Standard Hodge theory \cite{AmroucheEtAl1998,ArnoldFalkWinther2006} gives
\[
\dim\mathcal H_i=\beta_1(D_i),
\]
where \(\beta_1(D_i)\) denotes the first Betti number of \(D_i\).
Set \(N=\sum_i\beta_1(D_i)\). The constraint defining \(\mathcal Z_i\) is used only to construct local trial fields. The global eigenfields are not required to have zero normal trace on the inclusion interfaces.

Suppose that the inclusion tensor is \(A^\epsilon=\operatorname{diag}(1,1,\epsilon)\), so that \(e_z\) is the weak material direction. We start with local trial fields \(v_k=\phi_k(y)e_x\), where the \(\phi_k\) are smooth functions of one variable. Their curls satisfy \(\operatorname{curl}v_k=-\phi_k'(y)e_z\). 
To enforce the local orthogonality condition, find \(p_k\in H^1(D_i)/\mathbb R\) such that
\[
(\nabla p_k,\nabla q)_{D_i}=(v_k,\nabla q)_{D_i}
\qquad\text{for all }q\in H^1(D_i).
\]
The corrected field \(v_k-\nabla p_k\) then belongs to \(\mathcal Z_i\) and has the same curl as \(v_k\).
Subtracting the \(L^2(D_i)\)-orthogonal projection of each corrected field onto \(\mathcal H_i\) also preserves its curl. By choosing smooth functions \(\phi_k\) whose derivatives are linearly independent over the range of \(y\) in \(D_i\), we obtain weak trial spaces of arbitrarily large finite dimension in \(\mathcal Z_i\) that are orthogonal to \(\mathcal H_i\).

Let \(\mathcal S_K\subset\bigoplus_i\mathcal Z_i\) be the direct sum of all local harmonic spaces and \(K\) such weak fields, so \(\dim\mathcal S_K=N+K\). Bounded \(H(\operatorname{curl})\) extension operators \cite{hiptmair2020review}, combined with fixed smooth cutoffs around the separated inclusions, define a linear map \(T:\mathcal S_K\to H_0(\operatorname{curl};\Omega)\) satisfying \((Tw)|_{D_i}=w_i\) for every \(i\). By boundedness of the extension operators and norm equivalence on the fixed finite-dimensional space \(\mathcal S_K\),
\[
\|\operatorname{curl}Tw\|_{L^2(\Omega_0)}^2+
\sum_i\|\operatorname{curl}w_i\|_{L^2(D_i)}^2
\le C_K\|w\|_D^2
\qquad\text{for all }w\in\mathcal S_K.
\]
The constant \(C_K\) may depend on \(K\) and the chosen trial space, but is independent of \(\eta\) and \(\epsilon\).

The extension \(Tw\) satisfies the homogeneous tangential boundary condition on \(\partial\Omega\), but need not satisfy the weighted orthogonality condition defining \(\mathcal W_\eta\).
Define \(P^\eta v=v-\nabla p^\eta\), where \(p^\eta\in H_0^1(\Omega)\) is the unique solution of
\[
m^\eta(\nabla p^\eta,\nabla q)=m^\eta(v,\nabla q)
\quad\text{for all }q\in H_0^1(\Omega).
\]
It maps into \(\mathcal W_\eta\) and preserves curl. Although it need not preserve the inclusion fields themselves, their local orthogonality gives
\[
\begin{aligned}
m^\eta(P^\eta Tw,P^\eta Tw)
&\ge\eta\sum_i\|w_i-\nabla p^\eta\|_{L^2(D_i)}^2\\
&=\eta\sum_i\bigl(\|w_i\|_{L^2(D_i)}^2+
\|\nabla p^\eta\|_{L^2(D_i)}^2\bigr)
\ge\eta\|w\|_D^2.
\end{aligned}
\]
Indeed, \(p^\eta|_{D_i}\in H^1(D_i)\), so the cross term \((w_i,\nabla p^\eta)_{D_i}\) vanishes by the definition of \(\mathcal Z_i\). Thus \(P^\eta T\) is injective on \(\mathcal S_K\). Since all inclusion curls have only a weak component,
\[
a^{\eta,\epsilon}(P^\eta Tw,P^\eta Tw)
\le C_K(1+\eta\epsilon)\|w\|_D^2.
\]
Applying the min--max principle to the \((N+K)\)-dimensional trial space \(P^\eta T\mathcal S_K\subset\mathcal W_\eta\) gives
\[
\lambda^{\eta,\epsilon}_{N+K}\le C_K(\epsilon+\eta^{-1}).
\]
For \(N\ge1\), repeating the construction using only \(\bigoplus_i\mathcal H_i\) gives
\[
\lambda^{\eta,\epsilon}_{N}\le C\eta^{-1},
\]
since the local fields have zero curl in every inclusion.

\section{Finite element truncation and mesh dependence}
\label{sec:fem-truncation}
For fixed material parameters, the finite element eigenvalue approximation is governed by the standard spectral approximation framework \cite{Boffi2010}.
Let \(V_h\) be a conforming finite element space on an interface-fitted mesh. For Maxwell problems, we formulate the discrete eigenproblem on the weighted discrete complement \(\mathcal W_{\eta,h}\) defined in Section~\ref{sec:numerical-maxwell}. Let \(\mathcal K_h\) be an \(N\)-dimensional space of local discrete kernel fields, and let \(\mathcal W_h\subset\mathcal K_h^{\perp_{L^2(D)}}\) be a space of local discrete weak fields. 
In all models, the local discrete weak fields satisfy \(P_{\mathrm s}\mathcal Lw_h=0\) in \(D\). For Maxwell problems, we additionally require \(L^2(D_i)\)-orthogonality to the gradients of functions in the associated local scalar finite element space on each inclusion \(D_i\).
Set \(J_h=\dim\mathcal W_h\) and
\[
b_h(w_h,v_h)=\sum_i\int_{D_i}
\langle\mathcal Lw_h,q_s\rangle\langle\mathcal Lv_h,q_s\rangle\,dx.
\]

Assume that \(b_h\) is positive definite on \(\mathcal W_h\). For \(J_h\ge1\), denote its eigenvalues with respect to the local \(L^2(D)\) inner product, counted with multiplicity, by
\[
0<\mu_{h,1}\le\cdots\le\mu_{h,J_h}.
\]
More explicitly, for each \(j=1,\ldots,J_h\), the reduced eigenproblem seeks \(\psi_{h,j}\in\mathcal W_h\setminus\{0\}\) satisfying
\[
b_h(\psi_{h,j},v_h)=\mu_{h,j}(\psi_{h,j},v_h)_D
\qquad\text{for all }v_h\in\mathcal W_h,
\]
where \((u,v)_D=\sum_i\int_{D_i}u\cdot v\,dx\).
These eigenvalues measure local weak-direction energy before multiplication by \(\epsilon\) and before global extension.

\begin{theorem}[Discrete weak-mode upper bounds]
\label{thm:discrete-weak-bounds}
Suppose there is a linear discrete lifting operator \(T_h:\mathcal K_h\oplus\mathcal W_h\to V_h\), possibly depending on the material parameters, such that, for \(z_h=k_h+w_h\),
\begin{equation}
\label{eq:discrete-lifting}
\begin{aligned}
m^\eta(T_hz_h,T_hz_h)&\ge c\eta\|z_h\|_D^2,\\
\|\mathcal LT_hz_h\|_{L^2(\Omega_0)}^2
&\le C\bigl(\|z_h\|_D^2+b_h(w_h,w_h)\bigr),\\
\sum_i\int_{D_i}\langle A^\epsilon\mathcal LT_hz_h,
\mathcal LT_hz_h\rangle\,dx
&\le C\epsilon b_h(w_h,w_h).
\end{aligned}
\end{equation}
Here \(c,C>0\) are independent of \(h,\eta,\epsilon\), and \(V_h\) denotes the relevant unconstrained or constrained discrete space. The first inequality prevents loss of local mass and makes the lifting injective. The second controls the background energy of the lifted field, and the third bounds the energy contribution from strongly penalized components within the inclusions. Then
\begin{equation}
\label{eq:discrete-kernel-bound}
\lambda^{\eta,\epsilon}_{h,N}\le C'\eta^{-1}\quad(N\ge1),
\end{equation}
and
\begin{equation}
\label{eq:discrete-weak-bound}
\lambda^{\eta,\epsilon}_{h,N+j}
\le C'\eta^{-1}+C'(\epsilon+\eta^{-1})\mu_{h,j},
\qquad 1\le j\le J_h.
\end{equation}
If, in addition, \(\mu_{h,1}\ge c_0>0\) uniformly in \(h\), the latter estimate simplifies to
\[
\lambda^{\eta,\epsilon}_{h,N+j}
\le C''(\epsilon+\eta^{-1})\mu_{h,j}.
\]
\end{theorem}
\begin{proof}
Let \(\psi_{h,1},\ldots,\psi_{h,j}\) be local \(L^2\)-orthonormal eigenfunctions of \(b_h\), and set
\[
S_{h,j}=T_h\bigl(\mathcal K_h\oplus
\operatorname{span}\{\psi_{h,1},\ldots,\psi_{h,j}\}\bigr).
\]
The mass estimate makes \(T_h\) injective. For every \(z_h=k_h+w_h\) in this local trial space, we have \(b_h(w_h,w_h)\le\mu_{h,j}\|w_h\|_D^2\). The \(L^2(D)\)-orthogonality of \(k_h\) and \(w_h\) also gives \(\|w_h\|_D^2\le\|z_h\|_D^2\).
Hence
\[
\sup_{0\ne v_h\in S_{h,j}}
\frac{a^{\eta,\epsilon}(v_h,v_h)}{m^\eta(v_h,v_h)}
\le C'\eta^{-1}+C'(\epsilon+\eta^{-1})\mu_{h,j}.
\]
Applying the min--max principle to \(S_{h,j}\) yields the bound for \(\lambda^{\eta,\epsilon}_{h,N+j}\). Taking \(w_h=0\) and applying the same principle to \(T_h\mathcal K_h\) gives the kernel estimate when \(N\ge1\). The uniform lower bound \(\mu_{h,j}\ge\mu_{h,1}\ge c_0\) allows the additive \(\eta^{-1}\) term to be absorbed into \((\epsilon+\eta^{-1})\mu_{h,j}\), yielding the simplified bound.
\end{proof}

For finite elements of fixed polynomial degree on a shape-regular, quasi-uniform family of meshes, the standard inverse estimate gives
\[
b_h(w_h,w_h)\le Ch^{-2m}\|w_h\|_D^2
\quad\text{for all }w_h\in\mathcal W_h,
\]
with \(C\) independent of \(h\). Consequently, when \(J_h\ge1\), the largest reduced eigenvalue satisfies
\[
\mu_{h,J_h}\le Ch^{-2m}.
\]
Then, for \(0<h\le1\) and \(J_h\ge1\),
\[
\lambda^{\eta,\epsilon}_{h,N+J_h}
\le C(\epsilon+\eta^{-1})h^{-2m}.
\]
\Cref{sec:numex} provides explicit dimension counts for the discrete weak spaces used in the numerical experiments.

\begin{remark}
The dimension of the discrete space satisfying \(P_{\mathrm s}\mathcal Lw_h=0\) may depend on the alignment of the mesh with the weak material directions. Interpolation on a nonaligned mesh need not preserve this condition and may introduce components in the strongly penalized material subspace.
\end{remark}

\section{Spectral decoupling and eigenvalue counting for anisotropic diffusion}
\label{sec:modified-weyl}
For separable scalar diffusion, we derive an exact local counting formula and its small-\(\epsilon\) asymptotics, and establish the global decoupled limit as \(\eta\to\infty\) at fixed \(\epsilon>0\).

\subsection{A separable local model}
For scalar diffusion, let \(B^N_{i,\epsilon}\) be the self-adjoint operator associated with
\[
b_{i,\epsilon}(u,v)=\int_{D_i}A_{i,\epsilon}\nabla u\cdot\nabla v\,dx,
\qquad u,v\in H^1(D_i).
\]
Its boundary condition is the natural Neumann condition.

\begin{assumption}[Separable weak direction]
\label[assumption]{ass:separable-weak-direction}
For each \(i\), assume that in suitable Cartesian coordinates \(x=(y,z)\),
\[
D_i=Y_i\times(0,L_i),
\]
where \(L_i>0\) and \(Y_i\subset\mathbb R^{d-1}\) is a bounded, connected Lipschitz domain. The weak direction is the \(z\)-direction, and the local diffusion tensor has the form
\[
A_{i,\epsilon}(y,z)=
\begin{pmatrix}
A_{i,\perp}(y)&0\\
0&\epsilon
\end{pmatrix},
\]
where \(A_{i,\perp}\) is symmetric, bounded, and uniformly elliptic, with bounds independent of \(\epsilon\).
\end{assumption}

Let $B_{i,\perp}^{N}$ be the transverse Neumann operator on $Y_i$,
with eigenvalue problem
\[
  \int_{Y_i} A_{i,\perp}\nabla_y\phi\cdot\nabla_y\psi\,dy =
  \mu \int_{Y_i}\phi\psi\,dy, \qquad \psi\in H^1(Y_i).
\]
Its eigenvalues, repeated according to multiplicity, are denoted by
\[
  0=\mu_{i,0}<\mu_{i,1}\leq\mu_{i,2}\leq\cdots.
\]
The weak-direction Neumann eigenvalues on $(0,L_i)$ are
\begin{equation}
  \nu_{i,k} = \left(\frac{k\pi}{L_i}\right)^2, \qquad k\in\mathbb N_0.
  \label{eq:axial-eigenvalues}
\end{equation}
We write
\[
  N_{i,\epsilon}^{N}(\Lambda) := \#\left\{ \lambda\in\operatorname{spec}(B_{i,\epsilon}^{N}): \lambda\leq\Lambda \right\},
\]
where eigenvalues are counted with multiplicity.

\begin{theorem}[Exact counting formula and asymptotics as \(\epsilon\to0\)]
\label{thm:anisotropic-weyl}
Under \Cref{ass:separable-weak-direction}, the local operator admits the tensor-sum representation
\begin{equation}
  B_{i,\epsilon}^{N} = B_{i,\perp}^{N}\otimes I + \epsilon I\otimes \left(-\frac{d^2}{dz^2}\right)_{N}.
  \label{eq:tensor-sum}
\end{equation}
The two terms act on the transverse variables \(y\) and the weak-direction variable \(z\), respectively. Here \(I\) denotes the identity on the corresponding factor space, and \(\left(-d^2/dz^2\right)_N\) is the Neumann Laplacian on \((0,L_i)\).
Consequently, the eigenvalues of \(B_{i,\epsilon}^{N}\), counted with multiplicity, are
\begin{equation}
\operatorname{spec}(B_{i,\epsilon}^{N}) = \left\{\mu_{i,j}+\epsilon\nu_{i,k}: j,k\in\mathbb N_0 \right\},
  \label{eq:separated-spectrum}
\end{equation}
The corresponding counting function is given exactly by
\begin{equation}
\label{eq:exact-counting-law}
N_{i,\epsilon}^{N}(\Lambda) = \sum_{\substack{j\ge0\\ \mu_{i,j}\le\Lambda}} \left[ 1+ \left\lfloor \frac{L_i}{\pi\sqrt{\epsilon}} \sqrt{\Lambda-\mu_{i,j}} \right\rfloor \right].
\end{equation}
For every fixed \(\Lambda>0\), we have
\begin{equation}
  N_{i,\epsilon}^{N}(\Lambda) = \frac{L_i}{\pi\sqrt{\epsilon}} \sum_{\mu_{i,j}<\Lambda} \sqrt{\Lambda-\mu_{i,j}} +\mathcal O(1) \qquad\text{as }\epsilon\to0.
  \label{eq:modified-weyl-asymptotic}
\end{equation}
\end{theorem}

\begin{proof}
Let \(\{\phi_{i,j}\}_{j\ge0}\) be an orthonormal eigenbasis of \(B_{i,\perp}^{N}\) in \(L^2(Y_i)\), and let \(\{\chi_{i,k}\}_{k\ge0}\) be an orthonormal cosine eigenbasis of the Neumann Laplacian on \((0,L_i)\). Their products form an orthonormal basis of \(L^2(D_i)\). Substitution into the variational eigenproblem shows that \(\phi_{i,j}(y)\chi_{i,k}(z)\) has eigenvalue \(\mu_{i,j}+\epsilon\nu_{i,k}\). This establishes \eqref{eq:tensor-sum} and \eqref{eq:separated-spectrum}.

For each \(j\) with \(\mu_{i,j}\le\Lambda\), the indices \(k\in\mathbb N_0\) satisfying \(\mu_{i,j}+\epsilon\nu_{i,k}\le\Lambda\) are precisely those for which
\[
0\le k\le
\frac{L_i}{\pi\sqrt{\epsilon}}
\sqrt{\Lambda-\mu_{i,j}}.
\]
Their number is
\[
1+\left\lfloor
\frac{L_i}{\pi\sqrt{\epsilon}}
\sqrt{\Lambda-\mu_{i,j}}
\right\rfloor.
\]
Summing over \(j\) gives \cref{eq:exact-counting-law}.

Finally, \(0\le1+\lfloor t\rfloor-t\le1\) for \(t\ge0\). Hence
\[
N_{i,\epsilon}^{N}(\Lambda)
=
\frac{L_i}{\pi\sqrt{\epsilon}}
\sum_{\mu_{i,j}<\Lambda}
\sqrt{\Lambda-\mu_{i,j}}
+R_{i,\epsilon}(\Lambda),
\qquad
0\le R_{i,\epsilon}(\Lambda)
\le N_{i,\perp}^{N}(\Lambda).
\]
For fixed \(\Lambda\), \(N_{i,\perp}^{N}(\Lambda)\) is independent of \(\epsilon\). Thus \(R_{i,\epsilon}(\Lambda)=\mathcal O(1)\) as \(\epsilon\to0\), which proves the asymptotic formula.
\end{proof}

\begin{remark}
For a single inclusion and a threshold $0<\Lambda<\mu_{i,1}$, only the transverse mode $j=0$ contributes, and the weak-direction indices satisfy
\[
0\le k\le \left\lfloor \frac{L_i}{\pi}\sqrt{\frac{\Lambda}{\epsilon}} \right\rfloor .
\]
Each fixed positive index $k$ gives a local eigenvalue of order $\epsilon$, corresponding to the weak fields used in Section~4. The growing count at fixed $\Lambda$ arises because the largest admitted index increases as $\epsilon$ decreases. For $\Lambda=C\epsilon$, this largest index remains bounded. More generally, for a shrinking threshold $\Lambda_\epsilon\to0$, the count diverges precisely when $\Lambda_\epsilon/\epsilon\to\infty$.
\end{remark}

\subsection{High-contrast spectral decoupling}
We next relate the local counting formula to the global spectrum in the high-contrast limit. Define the global eigenvalue counting function by
\[
N_{\eta,\epsilon}(\Lambda)
:=
\#\{n\in\mathbb N:\lambda_n^{\eta,\epsilon}\le\Lambda\}.
\]
Let $B_{\Omega_0}^{D}$ denote the background operator on $\Omega_0$ with homogeneous Dirichlet conditions on both $\partial\Omega$ and the inclusion interfaces $\partial D_i$, and let $N_{\Omega_0}^{D}$ be its counting function.
Let \(B_{\eta,\epsilon}\) be the self-adjoint operator associated with the form \(a^{\eta,\epsilon}\) in the weighted space \(L^2(\Omega,\rho_\eta\,dx)\).

To compare the global operators as \(\eta\) varies, we use the fixed Hilbert space
\[
\mathcal H = L^2(\Omega_0)\oplus\bigoplus_{i=1}^M L^2(D_i).
\]

The map \(U_\eta\) defined below is unitary from \(L^2(\Omega,\rho_\eta\,dx)\) onto \(\mathcal H\). Define the transformed operator on \(\mathcal H\) by
\[
\widetilde B_{\eta,\epsilon} = U_\eta B_{\eta,\epsilon}U_\eta^{-1}.
\]

\begin{theorem}[High-contrast spectral decoupling]
\label{thm:high-contrast-decoupling}
Suppose that \Cref{ass:separable-weak-direction} holds,
and fix \(\epsilon>0\).
Under the unitary identification
\[
  U_\eta u
  :=
  \left(
    u|_{\Omega_0},
    \sqrt{\eta}\,u|_{D_1},
    \ldots,
    \sqrt{\eta}\,u|_{D_M}
  \right),
\]
the transformed operator
\(\widetilde B_{\eta,\epsilon}
=U_\eta B_{\eta,\epsilon}U_\eta^{-1}\)
converges in the norm-resolvent sense, as \(\eta\to\infty\), to
\begin{equation}
  B_{\infty,\epsilon}
  :=
  B_{\Omega_0}^{D}
  \oplus
  \bigoplus_{i=1}^{M}B_{i,\epsilon}^{N}.
  \label{eq:decoupled-limit-operator}
\end{equation}
Consequently, for every \(\Lambda>0\) satisfying
\(\Lambda\notin\operatorname{spec}(B_{\infty,\epsilon})\),
\begin{equation}
  \lim_{\eta\to\infty}N_{\eta,\epsilon}(\Lambda)
  =
  N_{\Omega_0}^{D}(\Lambda)
  +
  \sum_{i=1}^{M}N_{i,\epsilon}^{N}(\Lambda).
  \label{eq:global-counting-decoupling}
\end{equation}
\end{theorem}

\begin{proof}
\emph{Step 1: Transformed form and interface condition.}
The map \(U_\eta\) transforms the weighted mass inner product into the standard inner product on \(\mathcal H\). Writing \(v=(v_0,v_1,\ldots,v_M)=U_\eta u\), the transformed energy is
\begin{equation}
  \widetilde a_{\eta,\epsilon}(v,v) = \int_{\Omega_0}|\nabla v_0|^2\,dx + \sum_{i=1}^{M}b_{i,\epsilon}(v_i,v_i).
  \label{eq:transformed-energy}
\end{equation}
The traces of \(u\) from the background and inclusion sides agree, giving
\begin{equation}
  \operatorname{Tr}v_0 = \eta^{-1/2}\operatorname{Tr}v_i \qquad\text{on }\partial D_i.
  \label{eq:scaled-transmission}
\end{equation}
Thus, the transformed form domain \(\mathcal V_\eta=U_\eta H_0^1(\Omega)\) is
\begin{equation}
\begin{aligned}
\mathcal V_\eta=\bigl\{\, & v\in H^1(\Omega_0)\oplus \textstyle\bigoplus_{i=1}^{M}H^1(D_i): \\
& \operatorname{Tr}_{\partial\Omega}v_0=0,\quad \operatorname{Tr}_{\partial D_i}v_0 =\eta^{-1/2}\operatorname{Tr}_{\partial D_i}v_i \ \text{for all }i\,\bigr\}.
\end{aligned}
\label{eq:transformed-form-domain}
\end{equation}
For fixed \(\epsilon>0\), the energy expression contains no explicit dependence on \(\eta\). The transformed form depends on \(\eta\) only through the interface condition defining \(\mathcal V_\eta\).

\emph{Step 2: Uniform bounds and limiting boundary conditions.}
We first show how the interface condition determines the boundary
conditions of any limit. Both \(\widetilde B_{\eta,\epsilon}\) and \(B_{\infty,\epsilon}\) are nonnegative and self-adjoint, so their shifted resolvents are well defined. Set
\[
  R_{\eta,\epsilon} :=(\widetilde B_{\eta,\epsilon}+I)^{-1}, \qquad R_{\infty,\epsilon} :=(B_{\infty,\epsilon}+I)^{-1}.
\]
Fix \(f\in\mathcal H\) and let \(v^\eta=R_{\eta,\epsilon}f\).
Then
\[
  \widetilde a_{\eta,\epsilon}(v^\eta,w) +(v^\eta,w)_{\mathcal H} =(f,w)_{\mathcal H} \qquad\text{for all }w\in\mathcal V_\eta.
\]
Taking \(w=v^\eta\) and applying the Cauchy--Schwarz inequality gives
\[
  \|v^\eta\|_{\mathcal H}\le\|f\|_{\mathcal H}, \qquad \widetilde a_{\eta,\epsilon}(v^\eta,v^\eta) \le\|f\|_{\mathcal H}^{2}.
\]
For fixed \(\epsilon>0\), ellipticity and
\eqref{eq:transformed-energy} therefore imply
\[
  \|v^\eta\|_{ H^1(\Omega_0)\oplus\bigoplus_{i=1}^{M}H^1(D_i)} \le C_\epsilon\|f\|_{\mathcal H}, \qquad \eta\ge1,
\]
where \(C_\epsilon\) is independent of \(\eta\) and \(f\). Given any sequence \(\eta_n\to\infty\), the Rellich compactness theorem provides a subsequence, not relabelled, along which \(v^{\eta_n}\) converges weakly in the product \(H^1\) space and strongly in \(\mathcal H\) to some \(v=(v_0,v_1,\ldots,v_M)\).

The trace estimate and \eqref{eq:scaled-transmission} give
\[
  \|\operatorname{Tr}v_0^\eta\|_{H^{1/2}(\partial D_i)} \le C_\epsilon\eta^{-1/2}\|f\|_{\mathcal H} \longrightarrow0.
\]
Weak \(H^1\) convergence and continuity of the trace operator imply \(\operatorname{Tr}v_0=0\) on each \(\partial D_i\). The outer zero trace is also preserved.
Consequently, every such limit belongs to
\begin{equation}
  \mathcal V_\infty := H_0^1(\Omega_0) \oplus \bigoplus_{i=1}^{M}H^1(D_i).
  \label{eq:limiting-form-domain}
\end{equation}
Here \(H_0^1(\Omega_0)\) imposes zero trace on both \(\partial\Omega\) and all \(\partial D_i\). The space \(\mathcal V_\infty\) is the form domain of \(B_{\infty,\epsilon}\): the background has Dirichlet boundary conditions, while the inclusion forms with domains \(H^1(D_i)\) yield the natural Neumann operators.

\emph{Step 3: Identification of the limit equation.}
We now construct test functions satisfying the finite-\(\eta\) interface condition and use them to identify the limiting resolvent equation.
Since \(\Omega_0\) is a Lipschitz domain and the inclusion boundaries are separated from one another and from \(\partial\Omega\), for each \(i\) there exists a bounded linear trace lifting
\[
  R_i:H^{1/2}(\partial D_i)\longrightarrow H^1(\Omega_0)
\]
such that, for every \(g\in H^{1/2}(\partial D_i)\), \(R_i g\) has trace \(g\) on \(\partial D_i\) and zero trace on \(\partial\Omega\) and all other inclusion boundaries.

For any \(w=(w_0,w_1,\ldots,w_M)\in\mathcal V_\infty\), define
\(w^\eta=(w_0^\eta,w_1^\eta,\ldots,w_M^\eta)\) by
\begin{equation}
  w_i^\eta:=w_i \quad (i=1,\ldots,M),\qquad w_0^\eta := w_0+\eta^{-1/2} \sum_{i=1}^{M}R_i(\operatorname{Tr}w_i).
  \label{eq:approximating-test-functions}
\end{equation}
Then \(w^\eta\in\mathcal V_\eta\), and the boundedness of the liftings gives \(w^\eta\to w\) strongly in the product \(H^1\) space as \(\eta\to\infty\).
Taking \(w^{\eta_n}\) as the test function in the variational equation for \(v^{\eta_n}\), we pass to the limit using the weak convergence of the solutions and the strong convergence of the test functions in the product \(H^1\) space. This gives
\[
  \int_{\Omega_0}\nabla v_0\cdot\nabla w_0\,dx +\sum_{i=1}^{M}b_{i,\epsilon}(v_i,w_i) +(v,w)_{\mathcal H} =(f,w)_{\mathcal H} \qquad\text{for all }w\in\mathcal V_\infty.
\]
This is the variational equation for \(v=R_{\infty,\epsilon}f\), whose solution is unique. Since every sequence \(\eta_n\to\infty\) admits a subsequence along which \(R_{\eta_n,\epsilon}f\) converges strongly in \(\mathcal H\) to this same limit, the whole family satisfies
\[
  R_{\eta,\epsilon}f \longrightarrow R_{\infty,\epsilon}f \quad\text{strongly in }\mathcal H \quad\text{as }\eta\to\infty \qquad\text{for every }f\in\mathcal H.
\]

\emph{Step 4: Norm-resolvent convergence and eigenvalue counting.}
The same energy estimate holds for \(R_{\infty,\epsilon}\). Thus both resolvents map the unit ball of \(\mathcal H\) into a common bounded subset of the product \(H^1\) space. By the Rellich compactness theorem, this subset is relatively compact in \(\mathcal H\). We combine this compactness with self-adjointness to obtain convergence in operator norm.

Set \(T_\eta:=R_{\eta,\epsilon}-R_{\infty,\epsilon}\). If \(\|T_\eta\|_{\mathcal H\to\mathcal H}\) did not tend to zero as \(\eta\to\infty\), there would exist \(\eta_n\to\infty\), unit vectors \(f_n\in\mathcal H\), and \(\delta>0\) such that
\[
  \|T_{\eta_n}f_n\|_{\mathcal H}\ge\delta \qquad\text{for all }n.
\]
The uniform bound in the product \(H^1\) space and its compact embedding into \(\mathcal H\) yield a subsequence, not relabeled, such that \(T_{\eta_n}f_n\to g\) strongly in \(\mathcal H\). For every fixed \(h\in\mathcal H\), self-adjointness and Step~3 give
\[
  (g,h)_{\mathcal H} = \lim_{n\to\infty}(f_n,T_{\eta_n}h)_{\mathcal H} =0.
\]
Hence \(g=0\), contradicting \(\|g\|_{\mathcal H}\ge\delta\).
Therefore,
\[
  \|(\widetilde B_{\eta,\epsilon}+I)^{-1} -(B_{\infty,\epsilon}+I)^{-1}\|_{\mathcal H\to\mathcal H} \longrightarrow0.
\]

Finally, fix \(\Lambda>0\) outside \(\operatorname{spec}(B_{\infty,\epsilon})\), and choose \(a<0\). This negative lower endpoint ensures that the zero eigenvalues of the limiting Neumann components are included in \((a,\Lambda)\).
Norm-resolvent convergence implies norm convergence of the corresponding spectral projections for \((a,\Lambda)\).
These projections, therefore, have the same finite rank for all sufficiently large \(\eta\) \cite{Kato1995}.
Moreover, \(\Lambda\notin \operatorname{spec}(\widetilde B_{\eta,\epsilon})\) for all sufficiently large \(\eta\).
Since the operators are nonnegative, these ranks are precisely their eigenvalue counts up to \(\Lambda\).
The direct-sum representation \eqref{eq:decoupled-limit-operator} then gives \eqref{eq:global-counting-decoupling}.
\end{proof}

\section{Numerical examples}
\label{sec:numex}
The experiments use constant material directions and interface-fitted Cartesian meshes.
Section~\ref{subsec:weak_spec_branch} examines low-energy spectra and mesh-dependent spectral transitions in four models.
Section~\ref{subsec:numerical_decoupling_weyl} examines spectral decoupling and eigenvalue counting for scalar diffusion.

\subsection{Low-energy spectra and mesh dependence}
\label{subsec:weak_spec_branch}

\begin{figure}[htbp]
 \centering
 \includegraphics[width=0.48\linewidth]{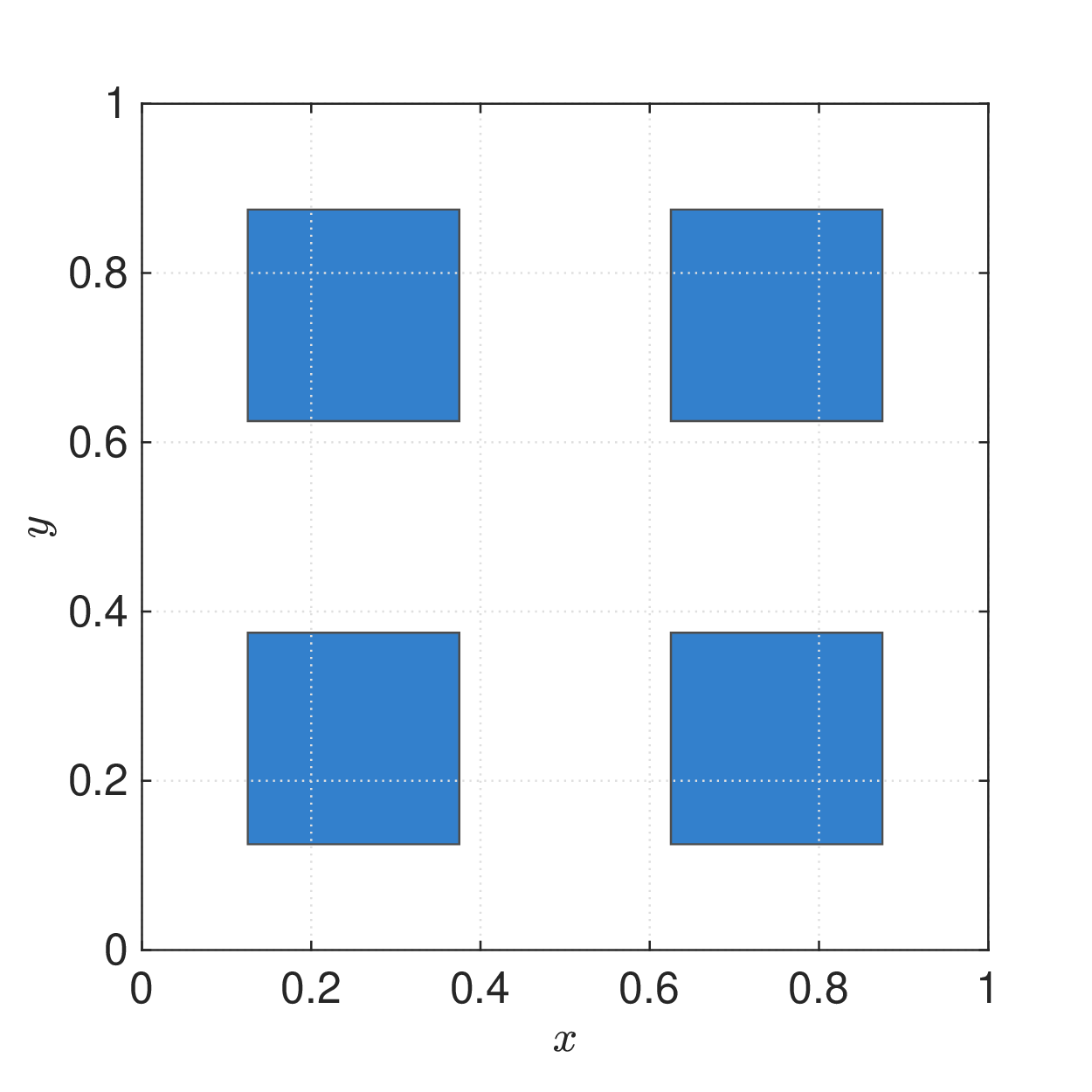}
 \caption{Geometry shared by the diffusion, elasticity, and biharmonic
experiments. The four shaded squares form the inclusion region $D$.}
 \label{fig:num-diffusion-geometry}
\end{figure}

\begin{table}[htbp]
\centering
\renewcommand{\arraystretch}{1.6}
\caption{Discretizations and numbers of displayed eigenvalues in the weak-mode experiments.}
\label{tab:numerical-settings}
\begin{tabularx}{\textwidth}{|l|X|X|r|}
\hline
\textbf{Model} &
\textbf{Elements} &
\textbf{Cells per coordinate direction} $n$ &
\textbf{Number displayed} \\ \hline
Diffusion & Scalar $Q_1$ & 16, 32, 64, 128 & 120 \\ \hline
Elasticity & Vector $Q_1$ & 16, 32, 64, 128 & 240 \\ \hline
Biharmonic & Bogner--Fox--Schmit & 16, 32, 64, 128 & 240 \\ \hline
Maxwell & Lowest-order first-family N\'ed\'elec  \cite{Nedelec1980} & 12, 18, 24 & 160 positive \\ \hline
\end{tabularx}
\end{table}

In all four experiments, \(\eta=10^6\) and \(\epsilon\in\{10^{-3},10^{-6}\}\). The diffusion, elasticity, and biharmonic models use the same geometry with four square inclusions, as shown in \Cref{fig:num-diffusion-geometry}. The discretizations and the numbers of displayed eigenvalues are summarized in \Cref{tab:numerical-settings}.

\subsubsection{Anisotropic diffusion eigenproblem}
\label{subsec:num-diffusion}
We first examine the weak spectral branch for scalar anisotropic diffusion. Let $\Omega=(0,1)^2$ and let $D$ be the union of four square inclusions:
\[
 D=\bigcup_{p,q=1}^{2}(I_p\times I_q), \qquad I_1=(1/8,3/8),\quad I_2=(5/8,7/8).
\]
The inclusions are strictly separated from one another and from $\partial\Omega$, and together occupy $25\%$ of the domain. The diffusion tensor is $I_2$ in the background and $\eta\operatorname{diag}(1,\epsilon)$ in each inclusion. The weak material direction is therefore the $y$-direction in all four inclusions. The homogeneous Dirichlet boundary condition excludes the global constant zero mode.

For \(\epsilon=10^{-3}\), the computed spectrum shows a pronounced separation between the fourth and fifth eigenvalues. On the finest mesh, \(\lambda_{h,4}\approx1.64\times10^{-4}\) and \(\lambda_{h,5}\approx1.58\times10^{-1}\). For \(\epsilon=10^{-6}\), the corresponding values are \(1.36\times10^{-4}\) and \(3.61\times10^{-4}\). The ratio \(\lambda_{h,5}/\lambda_{h,4}\) decreases from approximately \(965\) to \(2.66\). Thus stronger anisotropy substantially weakens the separation after the first four eigenvalues.

\begin{figure}[htbp]
    \centering
    \begin{subfigure}[htbp]{0.49\textwidth}
        \centering
        \includegraphics[width=\linewidth]{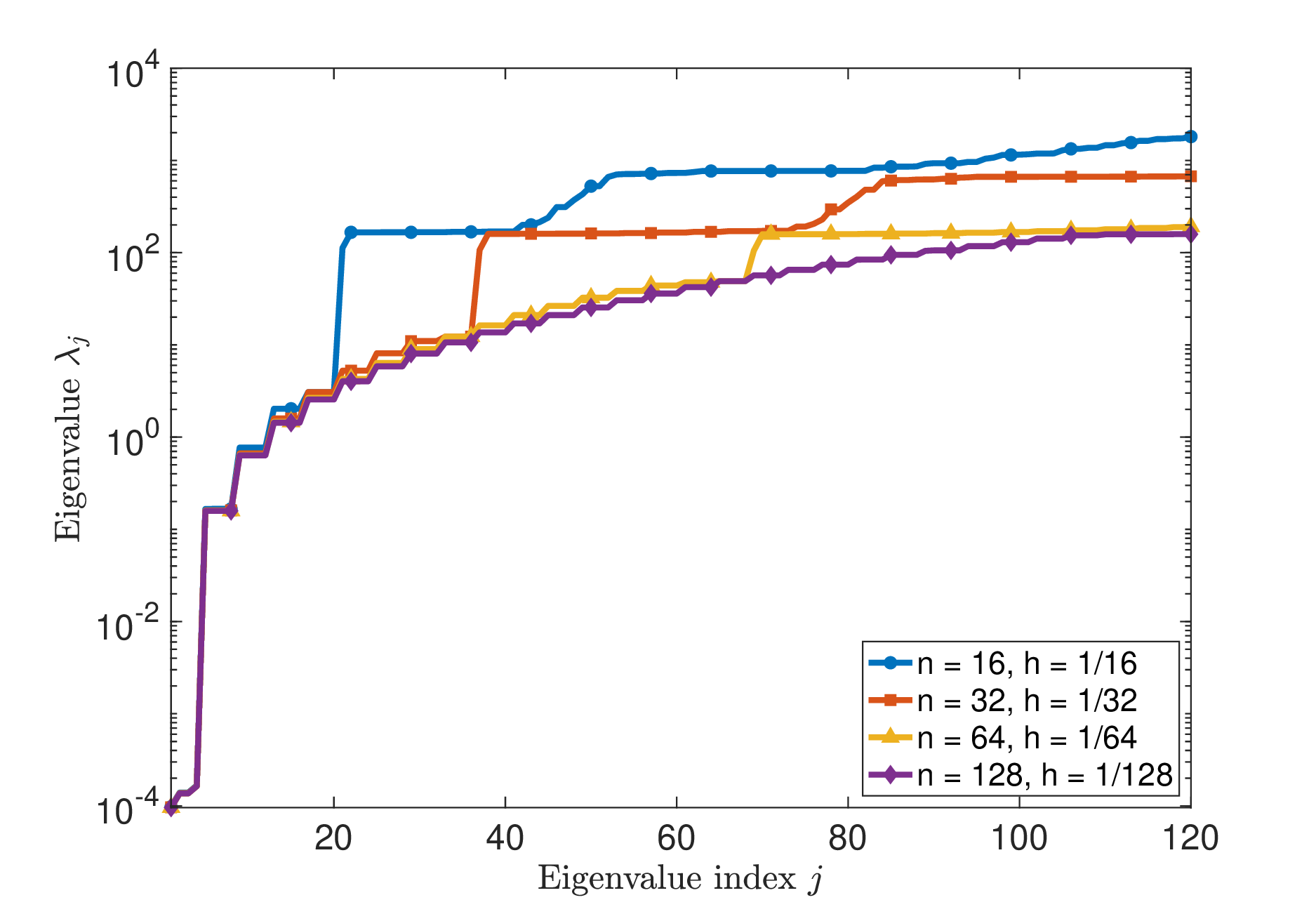}
        \caption{$\epsilon=10^{-3}$}
        \label{fig:num-diffusion-eps3}
    \end{subfigure}
    \hfill
    \begin{subfigure}[htbp]{0.49\textwidth}
        \centering
        \includegraphics[width=\linewidth]{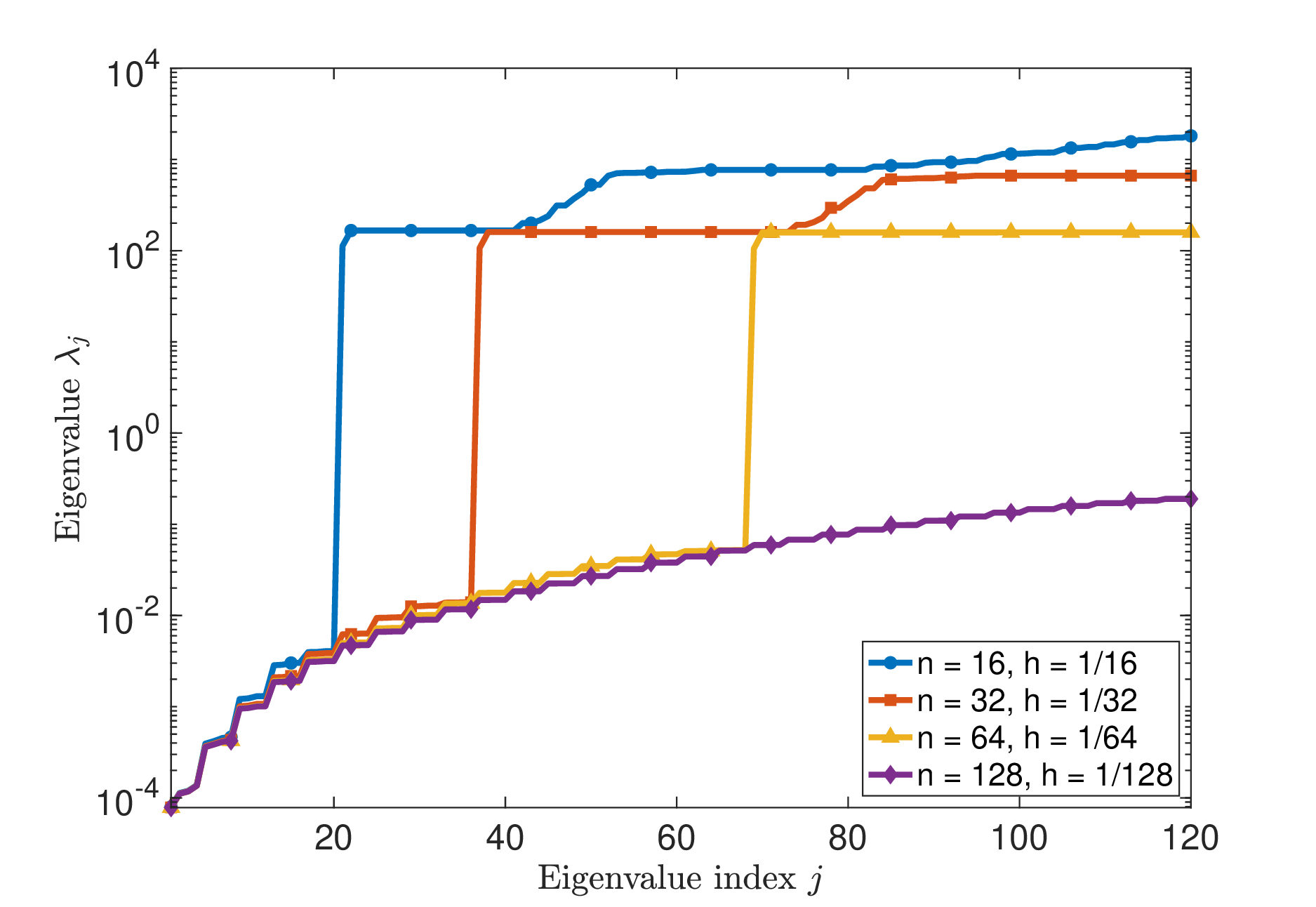}
        \caption{$\epsilon=10^{-6}$}
        \label{fig:num-diffusion-eps6}
    \end{subfigure}
    \caption{First $120$ eigenvalues of the anisotropic diffusion problem for $\eta=10^6$ on uniform meshes with $n=16,32,64,128$.}
    \label{fig:num-diffusion-spectra}
\end{figure}

The jumps in \Cref{fig:num-diffusion-eps6} occur after indices \(20,36,68\) for \(n=16,32,64\), respectively. Each inclusion contains \(n/4\) cells in the weak direction. The local \(Q_1\) fields depending only on this coordinate form a space of dimension \(n/4+1\). Orthogonality to the constants removes one dimension per inclusion, giving \(J_h=4[(n/4+1)-1]=n\). Since \(N=4\), the combined dimension is \(N+J_h=n+4\), which agrees with the observed jump locations. For \(n=128\), this dimension is \(132\), beyond the displayed range.

\subsubsection{Anisotropic elasticity eigenproblem}
\label{subsec:num-elasticity}
We next consider an elasticity problem with the same domain and four square inclusions as in Section~\ref{subsec:num-diffusion}. For a displacement field $\boldsymbol{u}=(u_1,u_2)^{\mathsf T}$, the Mandel strain vector is
\[
 \boldsymbol{e}_{\mathrm M}(\boldsymbol{u}) =
 \begin{pmatrix}
  \partial_x u_1\\
  \partial_y u_2\\
  (\partial_y u_1+\partial_x u_2)/\sqrt{2}
 \end{pmatrix}.
\]
We prescribe the material stiffness matrix in this representation as follows:
\begin{equation}
 C_{\eta,\epsilon}^{\mathrm M}(x) =
 \begin{cases}
  I_3, & x\in\Omega\setminus D,\\
  \eta\,\operatorname{diag}(1,1,\epsilon), & x\in D.
 \end{cases}
 \label{eq:num-elasticity-coefficient}
\end{equation}
Thus, only the shear strain component is weakly penalized inside the inclusions. The eigenproblem is to find $\lambda>0$ and
$0\ne\boldsymbol{u}\in[H_0^1(\Omega)]^2$ such that
\begin{equation}
 \int_\Omega \boldsymbol{e}_{\mathrm M}(\boldsymbol{v})^{\mathsf T} C_{\eta,\epsilon}^{\mathrm M} \boldsymbol{e}_{\mathrm M}(\boldsymbol{u})\,dx =
 \lambda\int_\Omega \rho_\eta\,\boldsymbol{u}\cdot\boldsymbol{v}\,dx,
 \label{eq:num-elasticity-eigenproblem}
\end{equation}
for all $\boldsymbol{v}\in[H_0^1(\Omega)]^2$.

In two dimensions, the rigid-motion space on each inclusion is
\[
 \mathcal{R} = \operatorname{span}
 \left\{ \begin{pmatrix}1\\0\end{pmatrix}, \begin{pmatrix}0\\1\end{pmatrix}, \begin{pmatrix}-y\\x\end{pmatrix} \right\}.
\]
Additional low-energy deformations arise from non-rigid fields of the form $\boldsymbol{u}=(\phi(y),\psi(x))^{\mathsf T}$ within an inclusion. Both normal strain components vanish, and the inclusion energy is
\[
 \frac{\eta\epsilon}{2} \int_{D_i}|\phi'(y)+\psi'(x)|^2\,dx\,dy.
\]

\begin{figure}[htbp]
 \centering
 \begin{minipage}[t]{0.49\linewidth}
  \centering
  \includegraphics[width=\linewidth] {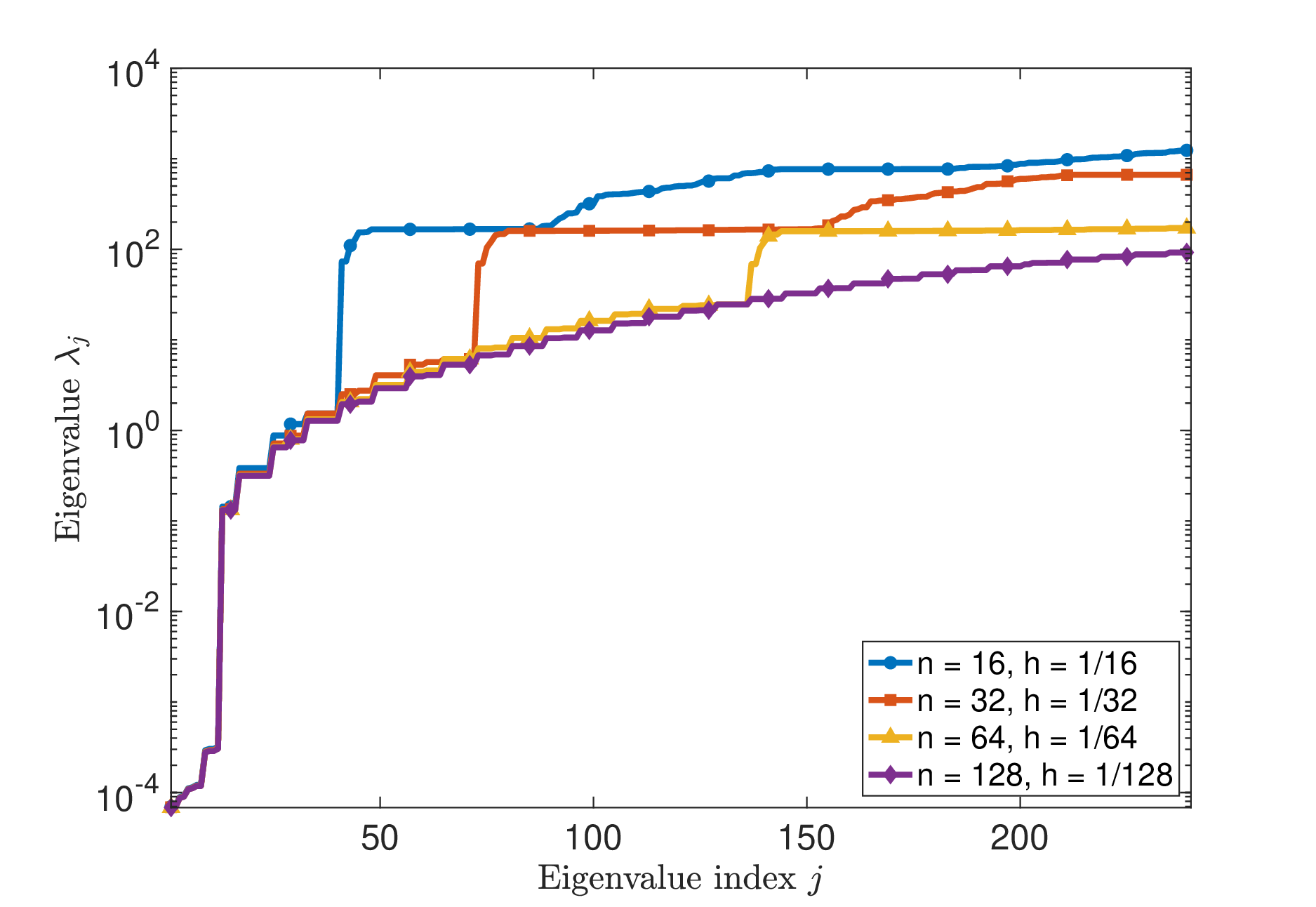}
  \par\smallskip
  {\small (a) $\epsilon=10^{-3}$}
 \end{minipage}\hfill
 \begin{minipage}[t]{0.49\linewidth}
  \centering
  \includegraphics[width=\linewidth] {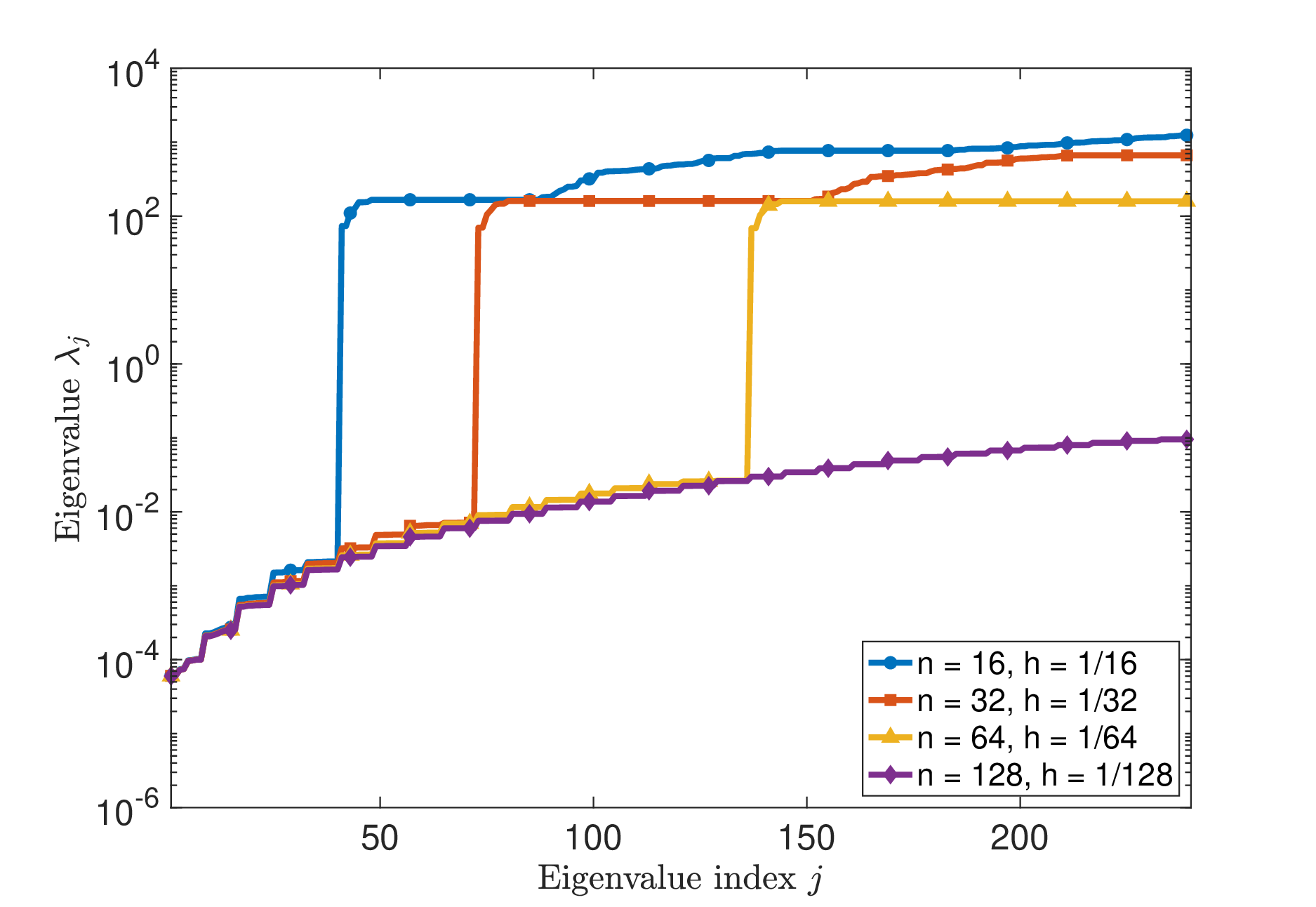}
  \par\smallskip
  {\small (b) $\epsilon=10^{-6}$}
 \end{minipage}
 \caption{First $240$ eigenvalues of the anisotropic elasticity problem for $\eta=10^6$ on uniform meshes with $n=16,32,64,128$.}
 \label{fig:num-elasticity-spectra}
\end{figure}

\Cref{fig:num-elasticity-spectra} shows that reducing \(\epsilon\) lowers the eigenvalues associated with weak shear deformations and weakens the separation after the first \(12\) eigenvalues. For \(\epsilon=10^{-6}\), the jumps occur after indices \(40,72,136\) for \(n=16,32,64\), respectively. With \(q=n/4\) cells per inclusion direction, the local \(Q_1\) displacement fields with zero normal strains form a space of dimension \(2(q+1)\). Removing the three rigid-motion directions on each inclusion gives
\[
J_h=4[2(q+1)-3]=2n-4,\qquad N+J_h=2n+8.
\]
The values of \(N+J_h\) match the observed jump locations. For \(n=128\), the corresponding index is \(264\), beyond the displayed range of the first \(240\) eigenvalues.

\subsubsection{Anisotropic biharmonic eigenproblem}
\label{subsec:num-biharmonic}
We use the same domain and four square inclusions as in Section~\ref{subsec:num-diffusion}.
In two dimensions, the Mandel representation of the Hessian is
\[
 \boldsymbol{H}_{\mathrm M}(u) = \bigl(u_{xx},\,u_{yy},\,\sqrt{2}\,u_{xy}\bigr)^{\mathsf T}.
\]
We prescribe the material coefficient matrix by
\[
 D_{\eta,\epsilon}^{\mathrm M}(x) =
 \begin{cases}
  I_3, & x\in\Omega\setminus D,\\
  \eta\,\operatorname{diag}(\epsilon,1,1), & x\in D.
 \end{cases}
\]
The eigenproblem is to find $\lambda>0$ and
$0\ne u\in H_0^2(\Omega)$ such that
\begin{equation}
 \int_\Omega \boldsymbol{H}_{\mathrm M}(v)^{\mathsf T} D_{\eta,\epsilon}^{\mathrm M} \boldsymbol{H}_{\mathrm M}(u)\,dx = \lambda\int_\Omega \rho_\eta uv\,dx,
 \label{eq:num-biharmonic-eigenproblem}
\end{equation}
for all $v\in H_0^2(\Omega)$. The outer boundary is therefore clamped: $u=\partial_n u=0$ on $\partial\Omega$. We use the conforming Bogner--Fox--Schmit element \cite{BrennerScott2008}.

The local Hessian kernel is the space of affine functions, \(\operatorname{span}\{1,x,y\}\), giving \(N=4\times3=12\) for the four inclusions.
Additional local weak fields are obtained by taking $u=\phi(x)$ inside an inclusion, with $\phi$ non-affine. For these fields, $u_{yy}=u_{xy}=0$, and the inclusion energy is
\[
 a_{D_i}(u,u) = \eta\epsilon\int_{D_i}|\phi''(x)|^2\,dx\,dy.
\]
Unlike the diffusion and elasticity energies, which involve first derivatives, the biharmonic energy involves second derivatives.
For an oscillatory function \(\phi\) with spatial frequency \(k\), the local Rayleigh quotient scales as \(\epsilon k^4\), compared with \(\epsilon k^2\) in the first-order models.

\begin{figure}[htbp]
 \centering
 \begin{minipage}[t]{0.49\linewidth}
  \centering
  \includegraphics[width=\linewidth]
   {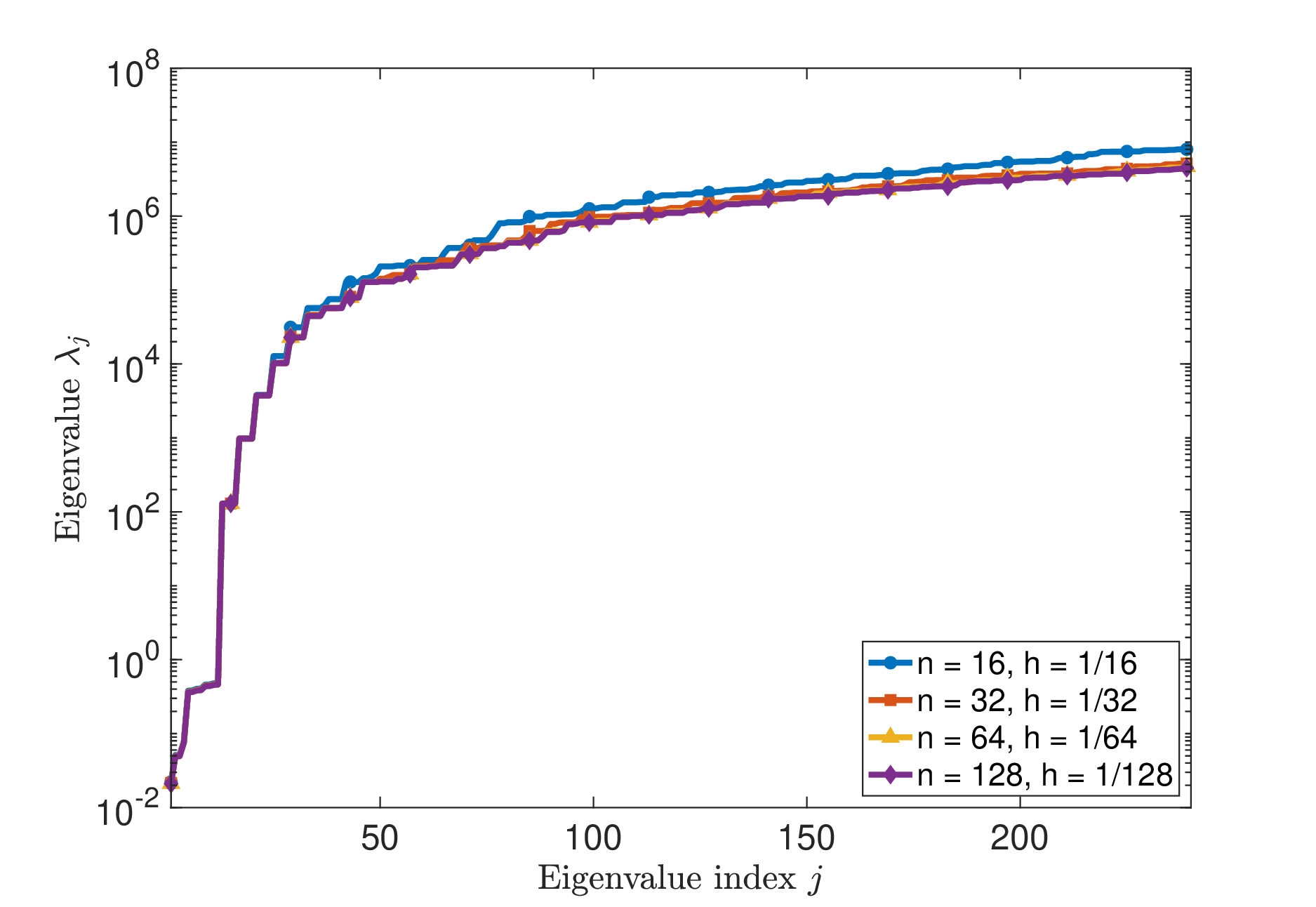}
  \par\smallskip
  {\small (a) $\epsilon=10^{-3}$}
 \end{minipage}\hfill
 \begin{minipage}[t]{0.49\linewidth}
  \centering
  \includegraphics[width=\linewidth]
   {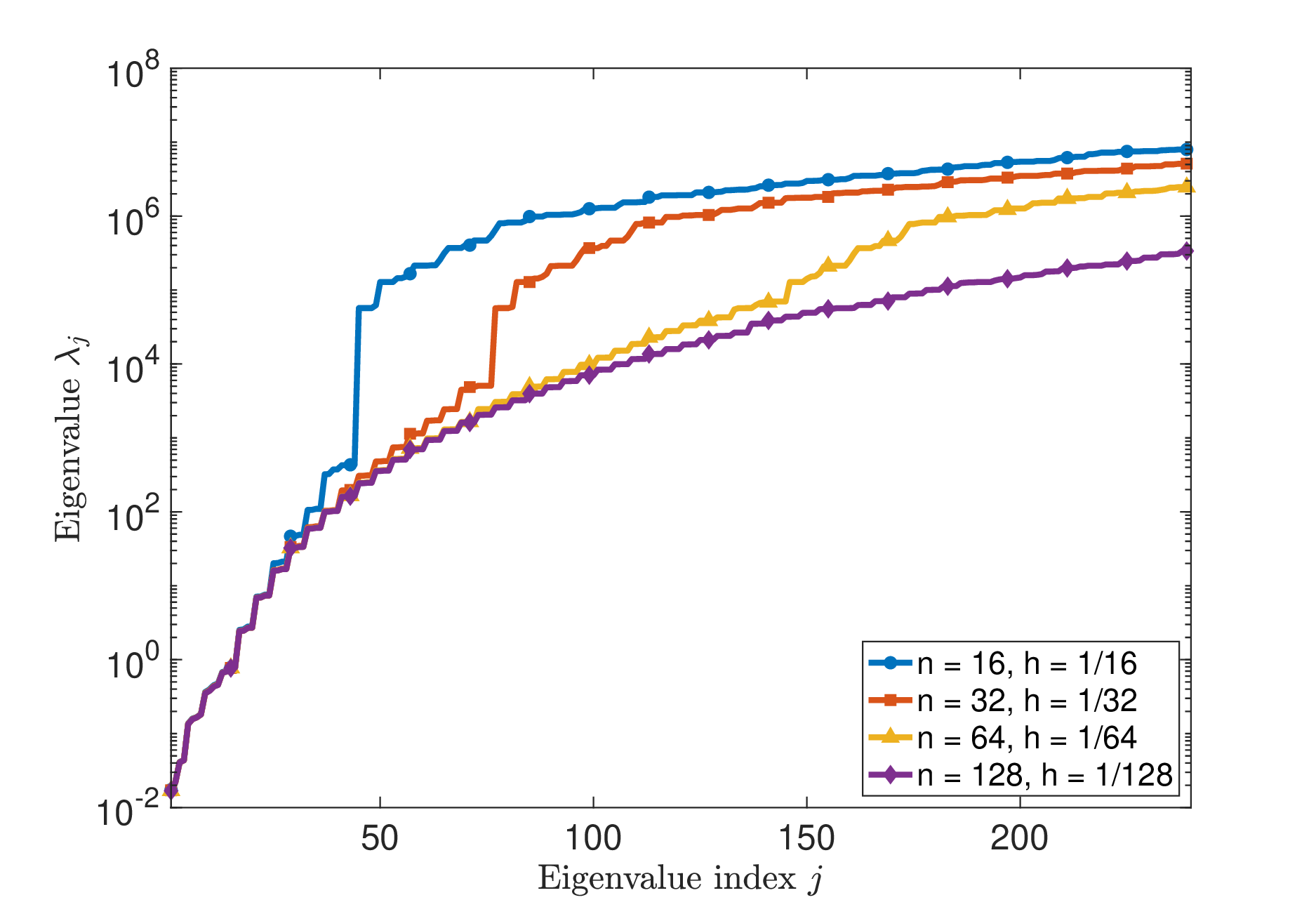}
  \par\smallskip
  {\small (b) $\epsilon=10^{-6}$}
 \end{minipage}
 \caption{First $240$ eigenvalues of the anisotropic biharmonic problem for $\eta=10^6$ on uniform meshes with $n=16,32,64,128$.}
 \label{fig:num-biharmonic-spectra}
\end{figure}

\Cref{fig:num-biharmonic-spectra} shows a downward shift of the low-energy spectrum and more pronounced mesh dependence as \(\epsilon\) decreases to \(10^{-6}\).
For the Bogner--Fox--Schmit space on the Cartesian meshes used here, the constraints \(u_{yy}=u_{xy}=0\) on each inclusion imply \(u_h=\phi_h(x)+cy\), where \(c\in\mathbb R\) and \(\phi_h\) is a one-dimensional \(C^1\) piecewise cubic function.
With \(q=n/4\) cells in the weak direction, this local constrained
space has dimension \(2(q+1)+1\).
Removing the three local kernel directions on each inclusion gives
\[
J_h=4[2(q+1)+1-3]=2n,\qquad N+J_h=2n+12.
\]
The \(\epsilon k^4\) scaling can raise the Rayleigh quotients of high-frequency weak fields to the background spectral range when \(\epsilon\) is not sufficiently small for the chosen mesh. This overlap can make the weak spectral branch less distinct.

\subsubsection{Anisotropic Maxwell eigenproblem}
\label{sec:numerical-maxwell}

\begin{figure}[htbp]
    \centering
    \includegraphics[width=0.48\textwidth]{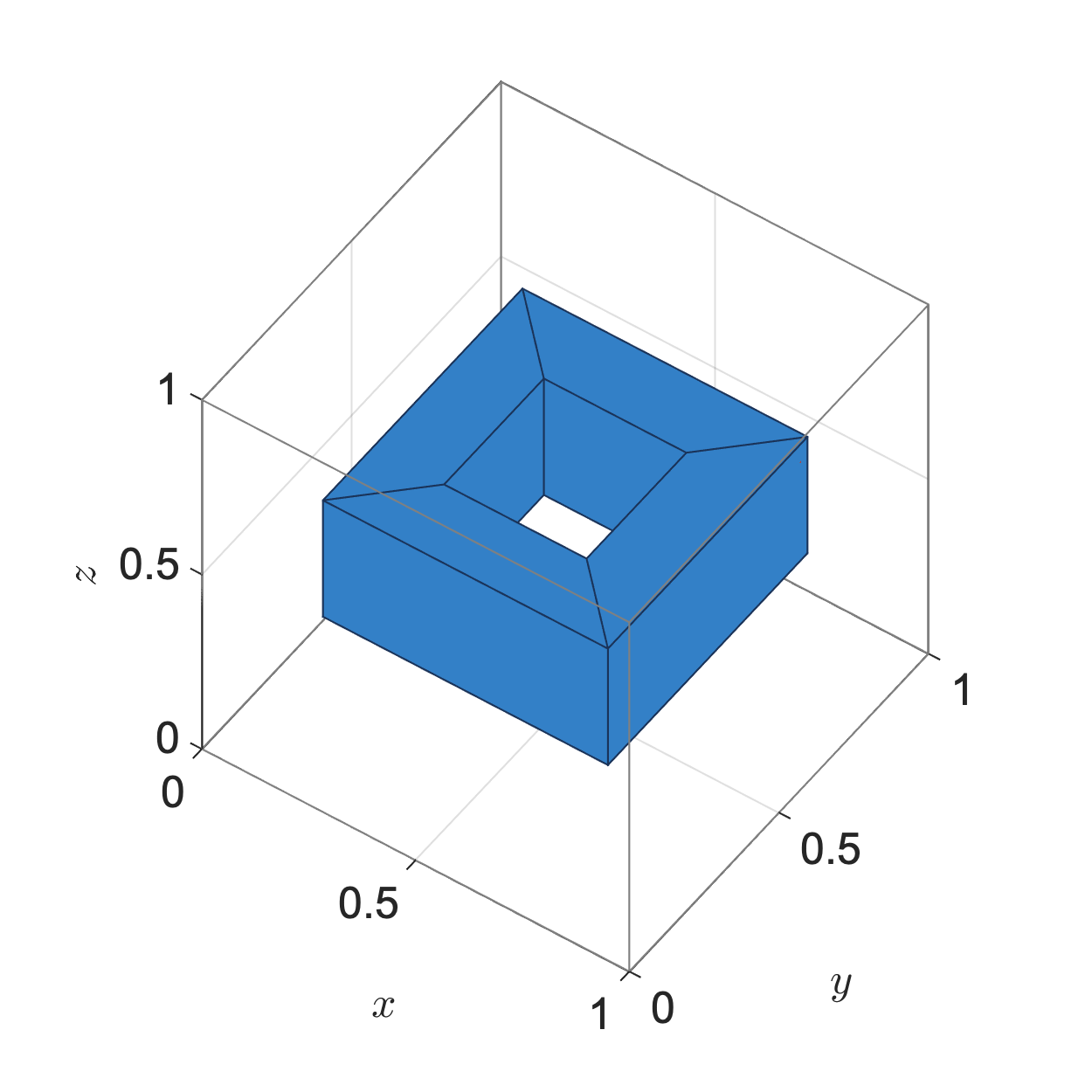}
    \caption{Geometry for the Maxwell experiment. The blue region represents the anisotropic inclusion $D$, and the rest of $\Omega$, including the central hole, is filled with background material.}
    \label{fig:maxwell-geometry}
\end{figure}

We next consider the three-dimensional Maxwell eigenvalue problem
on $\Omega=(0,1)^3$.
To obtain a multiply connected inclusion, we extrude a square
annulus in the $z$-direction and define
\[
D=\left(Q_{\mathrm{out}}\setminus\overline{Q_{\mathrm{in}}}\right) \times\left(\frac13,\frac23\right), \qquad
Q_{\mathrm{out}}=\left(\frac16,\frac56\right)^2,\quad
Q_{\mathrm{in}}=\left(\frac13,\frac23\right)^2.
\]
The geometry is shown in \Cref{fig:maxwell-geometry}. The inclusion has first Betti number $\beta_1(D)=1$, corresponding to one independent cohomology class associated with circulation around the hole. The hole is filled with background material.

The stiffness form is
\begin{align}
a^{\eta,\epsilon}(\boldsymbol{E},\boldsymbol{F}) &= \int_{\Omega\setminus D} (\nabla\times\boldsymbol{E})\cdot (\nabla\times\boldsymbol{F})\,dx \nonumber\\
&\quad+ \eta\int_D \bigl[\operatorname{diag}(1,1,\epsilon) (\nabla\times\boldsymbol{E})\bigr]\cdot (\nabla\times\boldsymbol{F})\,dx.
\label{eq:maxwell-numerical-stiffness}
\end{align}
Thus, the $z$-component of $\nabla\times\boldsymbol{E}$ is weakly penalized inside $D$. The perfectly conducting boundary condition $\boldsymbol{n}\times\boldsymbol{E}=0$ is imposed on $\partial\Omega$.

Let $V_h\subset H_0(\operatorname{curl};\Omega)$ denote the edge element space, and let $S_h\subset H_0^1(\Omega)$ be the associated continuous piecewise trilinear nodal space.
These spaces form part of a discrete de Rham complex, which provides the standard framework for treating the gradient kernel in Maxwell eigenvalue approximation \cite{Hiptmair2002,ArnoldFalkWinther2006}.
We remove the discrete gradient kernel by restricting the
eigenproblem to the weighted orthogonal complement
\[
\mathcal W_{\eta,h} = \left\{ \boldsymbol{v}_h\in V_h: m^\eta(\boldsymbol{v}_h,\nabla q_h)=0 \quad\text{for all }q_h\in S_h \right\}.
\]
We compute the positive eigenvalues $\lambda_h$ and the corresponding eigenfields $\boldsymbol{E}_h\in\mathcal W_{\eta,h}\setminus\{0\}$
that satisfy
\[
a^{\eta,\epsilon}(\boldsymbol{E}_h,\boldsymbol{F}_h) = \lambda_h m^\eta(\boldsymbol{E}_h,\boldsymbol{F}_h) \qquad \text{for all }\boldsymbol{F}_h\in\mathcal W_{\eta,h}.
\]

\begin{figure}[htbp]
    \centering
    \begin{subfigure}[t]{0.49\textwidth}
        \centering
        \includegraphics[width=\linewidth]{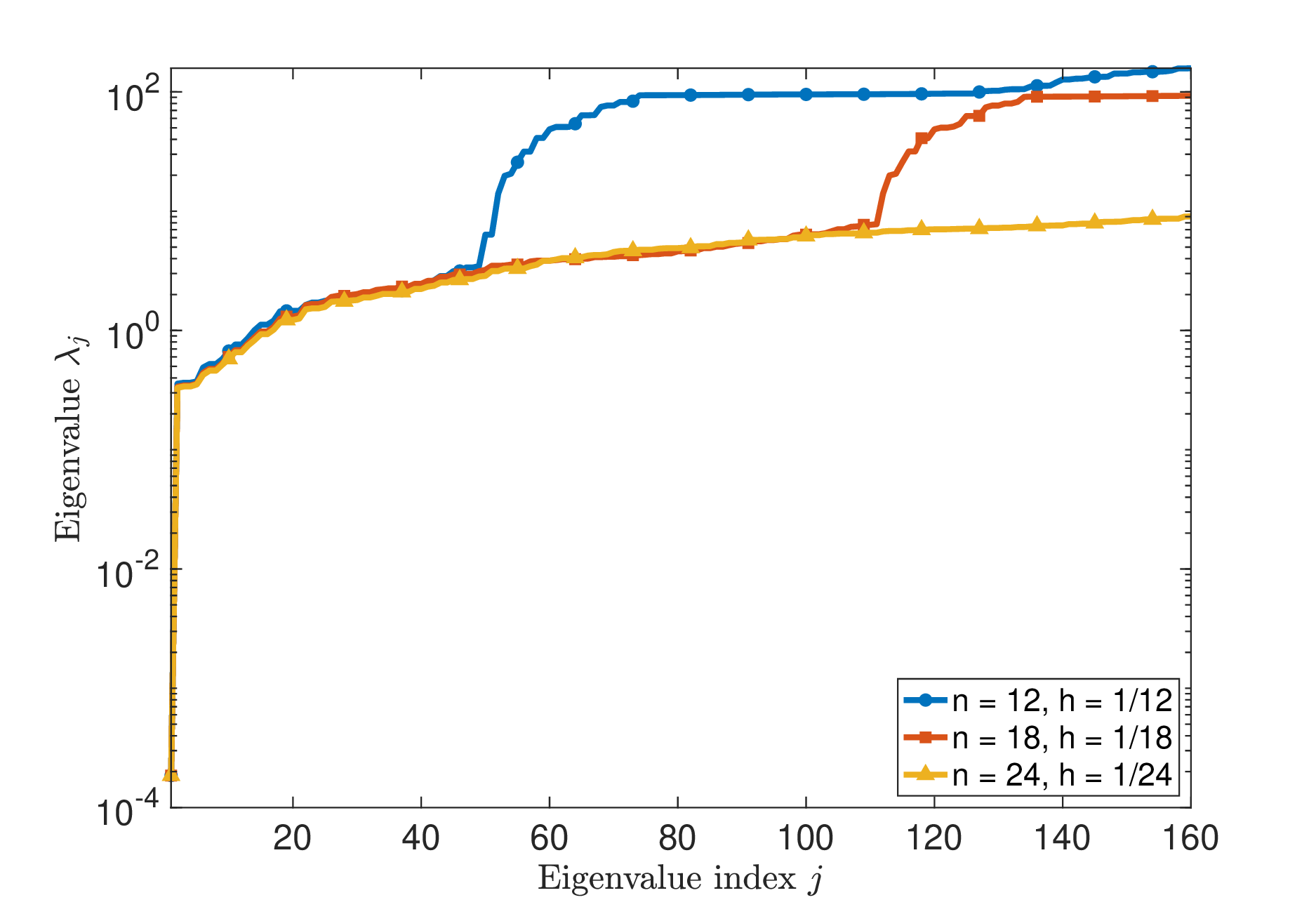}
        \caption{$\epsilon=10^{-3}$}
        \label{fig:maxwell-spectrum-eps3}
    \end{subfigure}
    \hfill
    \begin{subfigure}[t]{0.49\textwidth}
        \centering
        \includegraphics[width=\linewidth]{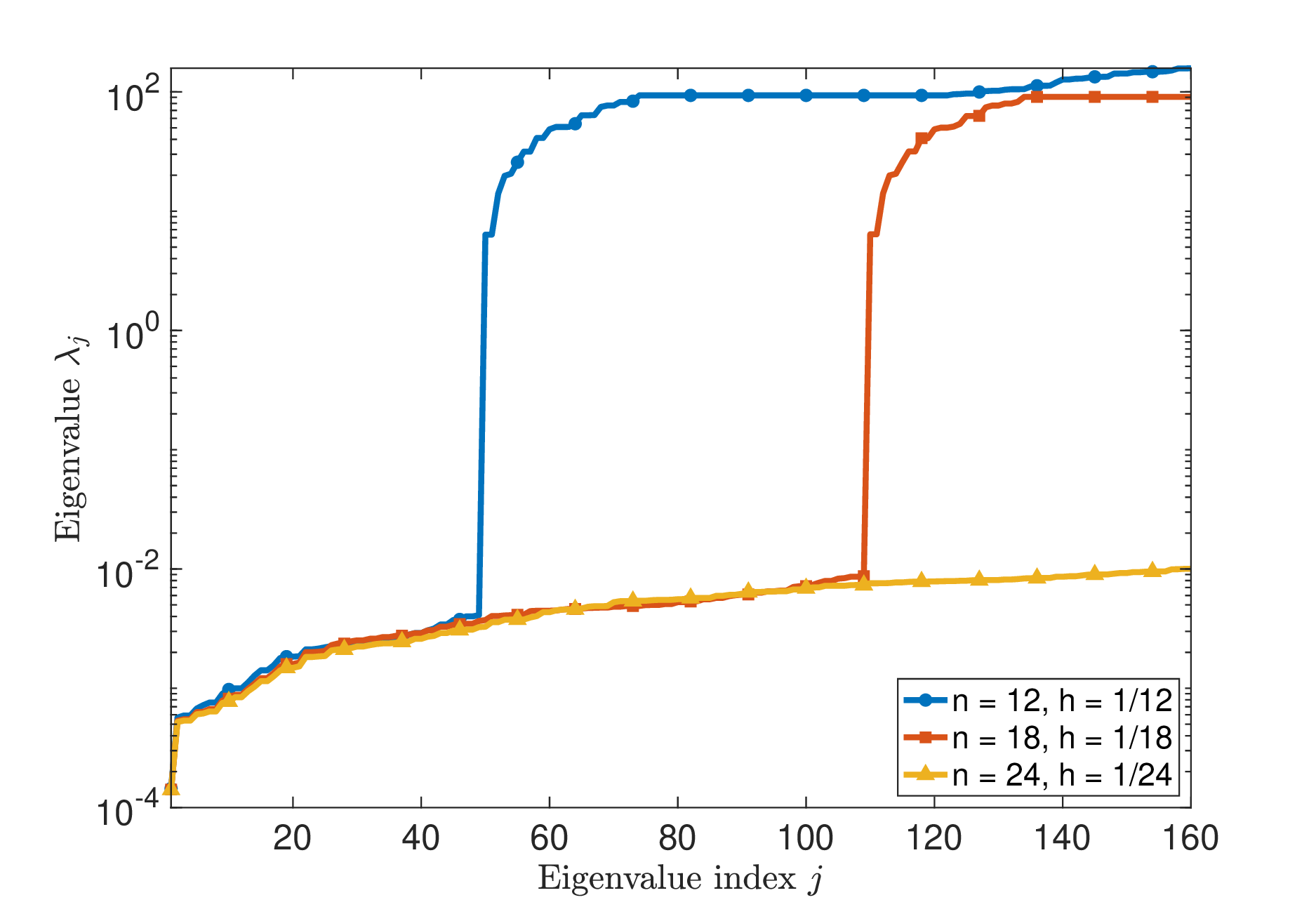}
        \caption{$\epsilon=10^{-6}$}
        \label{fig:maxwell-spectrum-eps6}
    \end{subfigure}
    \caption{First $160$ positive eigenvalues of the anisotropic Maxwell problem for $\eta=10^6$ on uniform meshes with $n=12,18,24$.}
    \label{fig:maxwell-spectra}
\end{figure}

\Cref{fig:maxwell-spectra} shows that the low-energy spectrum depends strongly on both anisotropy and mesh resolution. For $\epsilon=10^{-6}$, sharp spectral jumps occur after indices $49$ and $109$ for $n=12$ and $n=18$, respectively. These jump indices agree with a local dimension count that combines the weak curl degrees of freedom with the topological contribution of the inclusion.

Let
\[
\mathcal{R} = Q_{\mathrm{out}}\setminus\overline{Q_{\mathrm{in}}}
\]
denote the annular cross section of $D$. On the uniform interface-fitted Cartesian mesh with spacing $h=1/n$, the number of square cells in this cross section is
\[
F_{\mathcal{R}} = \left(\frac{2n}{3}\right)^2 - \left(\frac{n}{3}\right)^2 = 3\left(\frac{n}{3}\right)^2.
\]
Consider a local edge element field on $D$ of the form
\[
\boldsymbol{E}_h(x,y,z) = \bigl(E_{x,h}(x,y),E_{y,h}(x,y),0\bigr).
\]
Since this field is independent of $z$, its curl is
\[
\nabla\times\boldsymbol{E}_h = \left(0,0,\partial_x E_{y,h}-\partial_y E_{x,h}\right).
\]
Its curl is therefore aligned with the weak material direction inside $D$.

On the annular cross section $\mathcal R$, let $V_{\mathcal R},E_{\mathcal R},F_{\mathcal R}$ denote the numbers of vertices, edges, and cells, including boundary vertices and edges. No essential boundary conditions are imposed on the local nodal and edge element spaces. Since $\mathcal R$ is connected, the local discrete gradient has rank $V_{\mathcal R}-1$.
The Euler relation for the annular mesh is $V_{\mathcal R}-E_{\mathcal R}+F_{\mathcal R}=0$. Therefore, the planar edge space modulo discrete gradients has dimension
\[
N_h^{\mathrm{loc}} = E_{\mathcal R}-(V_{\mathcal R}-1) = F_{\mathcal R}+1 = 3(n/3)^2+1.
\]
The additional dimension corresponds to a local discrete harmonic field with nonzero circulation around the hole. For $n=12,18,24$, the resulting dimensions are $49,109,193$, respectively.
The first two values agree with the observed jump indices in \Cref{fig:maxwell-spectrum-eps6}. The third lies beyond the displayed range.

\subsection{High-contrast spectral decoupling and anisotropic eigenvalue counting}
\label{subsec:numerical_decoupling_weyl}
The scalar experiments test two predictions from \Cref{sec:modified-weyl}: decoupling as $\eta$ increases at fixed $\epsilon$, and eigenvalue counting in a fixed energy window as $\epsilon$ decreases.

\begin{figure}[htbp]
    \centering
    \includegraphics[width=0.44\textwidth]
    {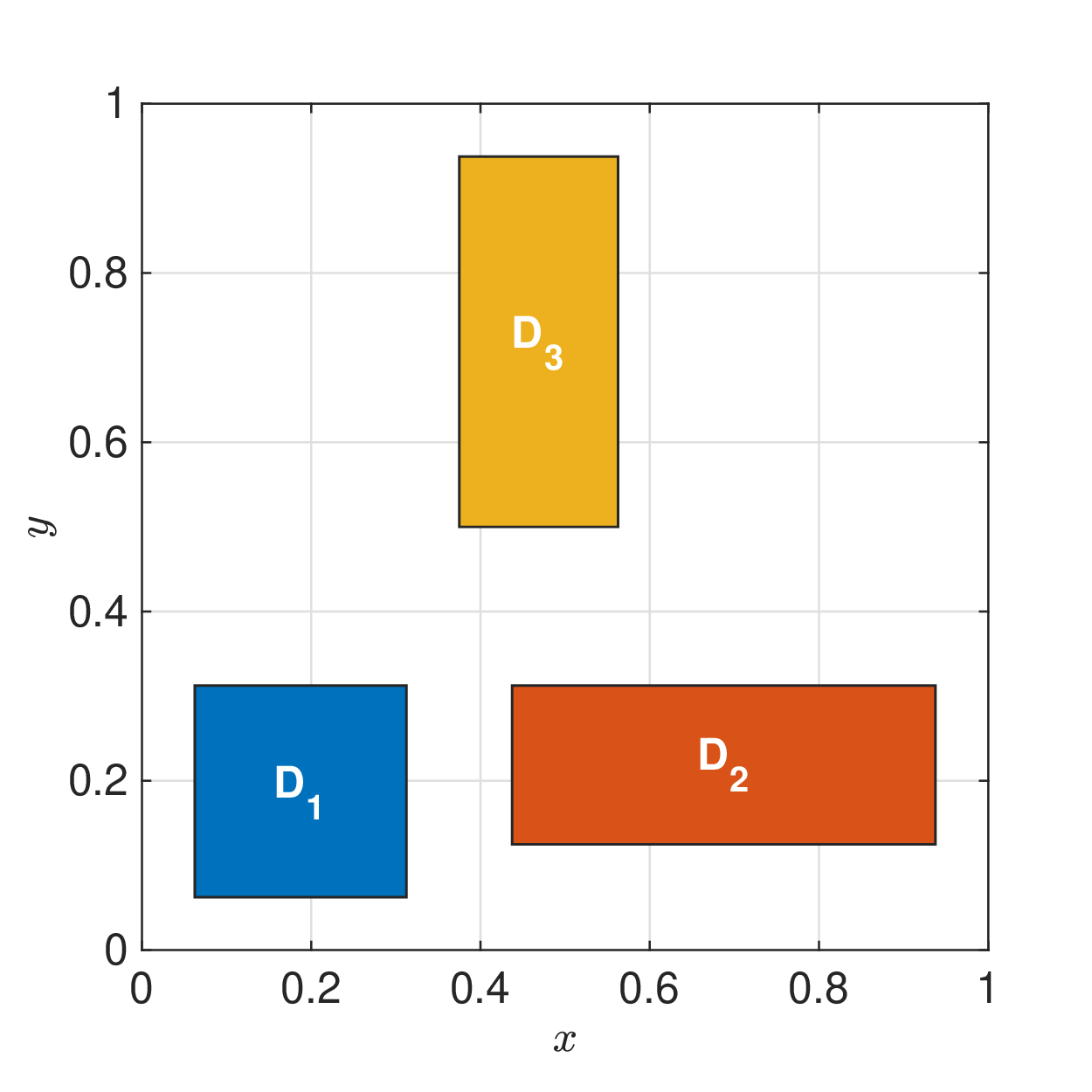}
    \caption{Three rectangular inclusions with different sizes and aspect ratios.}
    \label{fig:three-inclusion-geometry}
\end{figure}

Let $\Omega=(0,1)^2$ and consider three disjoint rectangular
inclusions with different sizes and aspect ratios:
\begin{align*}
D_1&=\left(\frac1{16},\frac5{16}\right) \times\left(\frac1{16},\frac5{16}\right),\\
D_2&=\left(\frac7{16},\frac{15}{16}\right) \times\left(\frac18,\frac5{16}\right),\\
D_3&=\left(\frac38,\frac9{16}\right) \times\left(\frac12,\frac{15}{16}\right).
\end{align*}
The geometry is illustrated in \Cref{fig:three-inclusion-geometry}.
Writing $D=D_1\cup D_2\cup D_3$ and
$\Omega_0=\Omega\setminus\overline{D}$, we set
\[
\kappa_{\eta,\epsilon}(x)=
\begin{cases}
\eta\operatorname{diag}(1,\epsilon), & x\in D,\\
I_2, & x\in\Omega_0.
\end{cases}
\]
The weak material direction is therefore the $y$-direction in all three inclusions.
We seek $\lambda>0$ and $u\in H_0^1(\Omega)\setminus\{0\}$ such that
\[
\int_\Omega \kappa_{\eta,\epsilon}\nabla u\cdot\nabla v\,dx = \lambda\int_\Omega \rho_\eta uv\,dx \qquad\text{for all }v\in H_0^1(\Omega).
\]
The discretization uses conforming bilinear finite elements on uniform Cartesian meshes that fit all inclusion interfaces.

We first fix $\epsilon=10^{-2}$ and vary
\[
\eta\in\{10^2,10^3,10^4,10^5,10^6\}.
\]
By \Cref{thm:high-contrast-decoupling}, the limiting operator
as $\eta\to\infty$ is
\[
B_{\infty,\epsilon} = B_{\Omega_0}^{D} \oplus \bigoplus_{i=1}^{3}B_{i,\epsilon}^{N}.
\]
Here, $B_{\Omega_0}^{D}$ has homogeneous Dirichlet conditions on $\partial\Omega$ and all inclusion interfaces, and $B_{i,\epsilon}^{N}$ has natural Neumann conditions on $\partial D_i$.

For each mesh, we construct a discrete reference spectrum by merging the finite element spectra of the background Dirichlet problem and the three local Neumann problems, retaining multiplicities. Denote these reference eigenvalues by $\lambda_{\epsilon,h,j}^{\mathrm{dec}}$ and the global finite element eigenvalues by $\lambda^{\eta,\epsilon}_{h,j}$. Both spectra use the same underlying mesh and finite element discretization, so their difference measures the finite-contrast coupling relative to the discrete decoupled limit. The first three reference eigenvalues are zero, corresponding to one constant Neumann mode in each inclusion.

\begin{figure}[htbp]
    \centering
    \begin{subfigure}[htbp]{0.49\textwidth}
        \centering
        \includegraphics[width=\linewidth] {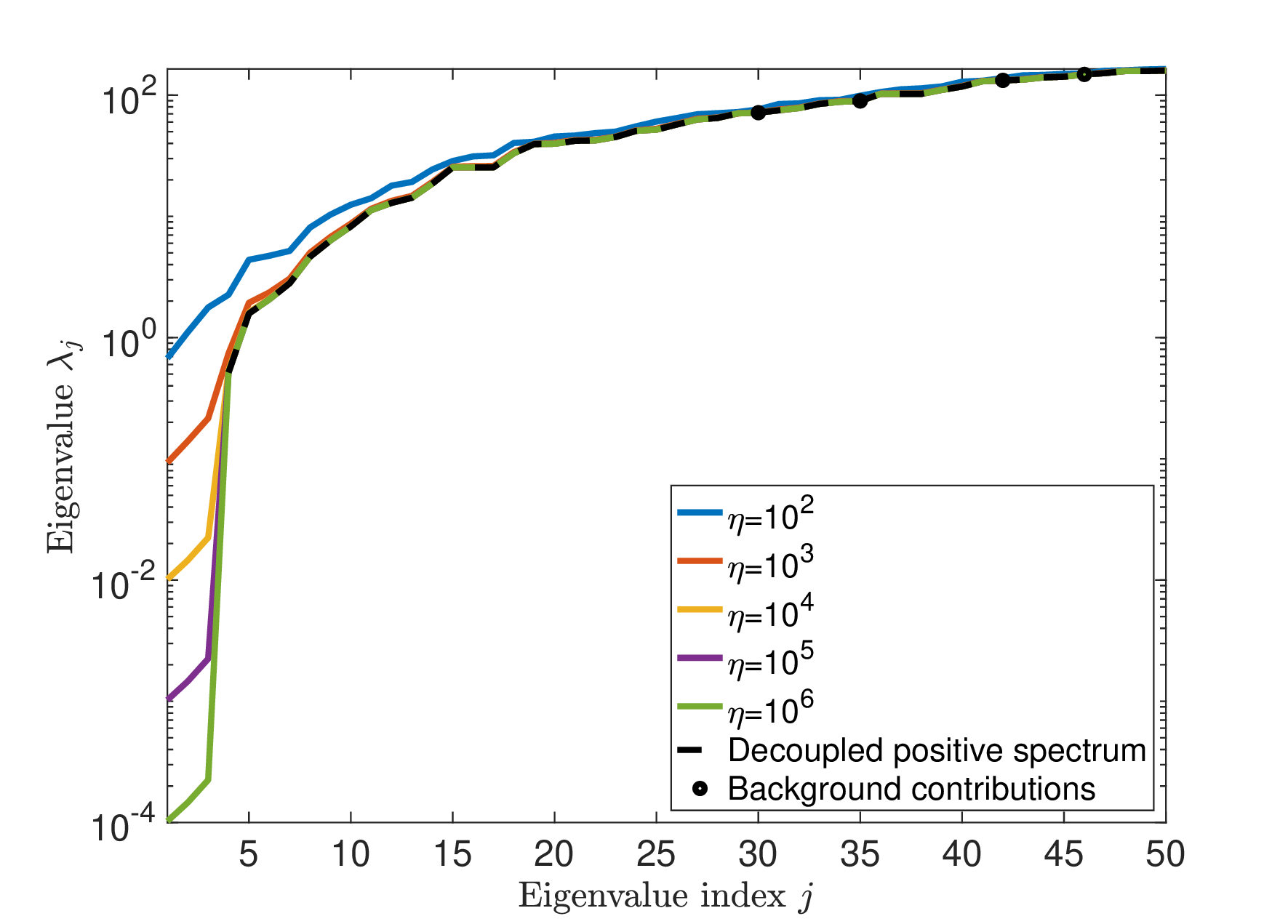}
        \caption{Global and decoupled spectra, $n=256$.}
        \label{fig:decoupling-spectrum}
    \end{subfigure}
    \hfill
    \begin{subfigure}[htbp]{0.49\textwidth}
        \centering
        \includegraphics[width=\linewidth]
        {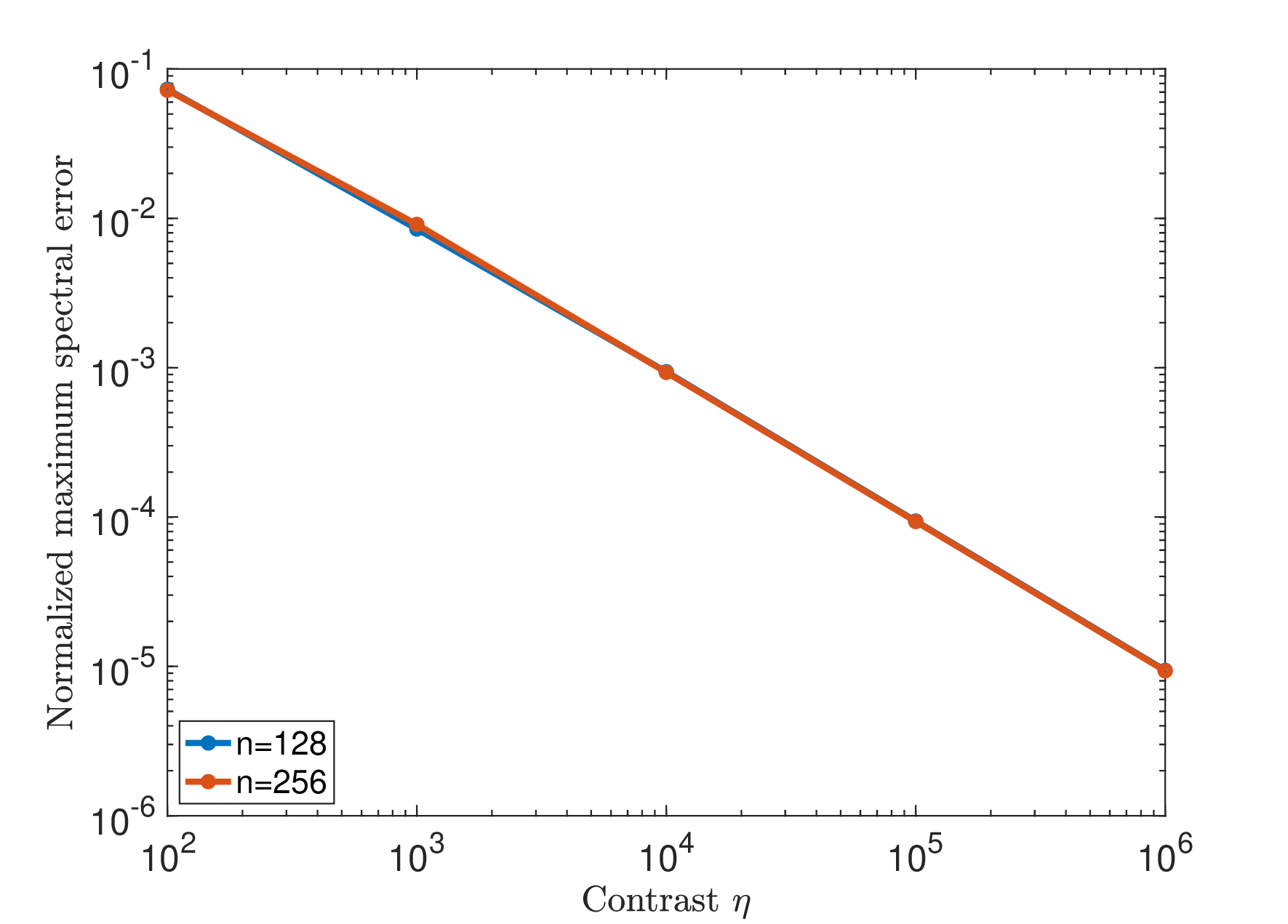}
        \caption{Normalized maximum spectral discrepancy.}
        \label{fig:decoupling-error}
    \end{subfigure}
    \caption{ High-contrast spectral decoupling for $\epsilon=10^{-2}$. Left: the first $50$ global eigenvalues for different values of $\eta$, compared with the positive part of the discrete decoupled reference spectrum. Right: the discrepancy defined in \eqref{eq:numerical-decoupling-error}.} \label{fig:decoupling-results}
\end{figure}

\Cref{fig:decoupling-spectrum} compares the first $50$ global eigenvalues with the positive part of the discrete decoupled reference spectrum on the mesh with $n=256$ cells in each coordinate direction. As $\eta$ increases, the first three global eigenvalues approach zero, while the remaining displayed eigenvalues approach their positive decoupled counterparts. The marked background Dirichlet eigenvalues illustrate that the decoupled spectrum includes contributions from both the inclusions and the background.

To quantify this convergence, we use the normalized maximum
spectral discrepancy
\begin{equation}
E_{\eta,h} = \frac{ \displaystyle\max_{1\le j\le K} \left| \lambda^{\eta,\epsilon}_{h,j} - \lambda_{\epsilon,h,j}^{\mathrm{dec}} \right|}{1+\lambda_{\epsilon,h,K}^{\mathrm{dec}}}, \qquad K=50.
\label{eq:numerical-decoupling-error}
\end{equation}
This comparison includes the three zero eigenvalues of the decoupled problem. \Cref{fig:decoupling-error} shows the results for $n=128$ and $n=256$. On the finer mesh, the discrepancy decreases from approximately $7.21\times10^{-2}$ at $\eta=10^2$ to $9.33\times10^{-6}$ at $\eta=10^6$.
For both meshes, $E_{\eta,h}$ decays approximately as $\eta^{-1}$ for large $\eta$. An $O(\eta^{-1})$ bound for $E_{\eta,h}$ can also be established in this setting. We omit the proof due to space limitations.

For comparison with the global finite element counts, we introduce
the counting function of the decoupled operator:
\[
N_{\infty,\epsilon}(\Lambda) = \#\left\{ \lambda\in\operatorname{spec}(B_{\infty,\epsilon}): \lambda\le\Lambda \right\},
\]
with eigenvalues counted with multiplicity.

Let $\ell_i$ and $L_i$ denote the side lengths of $D_i$ in the strong and weak directions, respectively. For the chosen geometry,
\[
(\ell_1,\ell_2,\ell_3) = \left(\frac14,\frac12,\frac3{16}\right), \qquad
(L_1,L_2,L_3) = \left(\frac14,\frac3{16},\frac7{16}\right).
\]
The local Neumann eigenvalues are
\begin{equation}
\lambda_{i,j,k}^{N} = \left(\frac{j\pi}{\ell_i}\right)^2 + \epsilon\left(\frac{k\pi}{L_i}\right)^2, \qquad j,k\in\mathbb N_0.
\label{eq:three-inclusion-local-spectrum}
\end{equation}

\begin{figure}[htbp]
    \centering
    \includegraphics[width=0.70\textwidth]{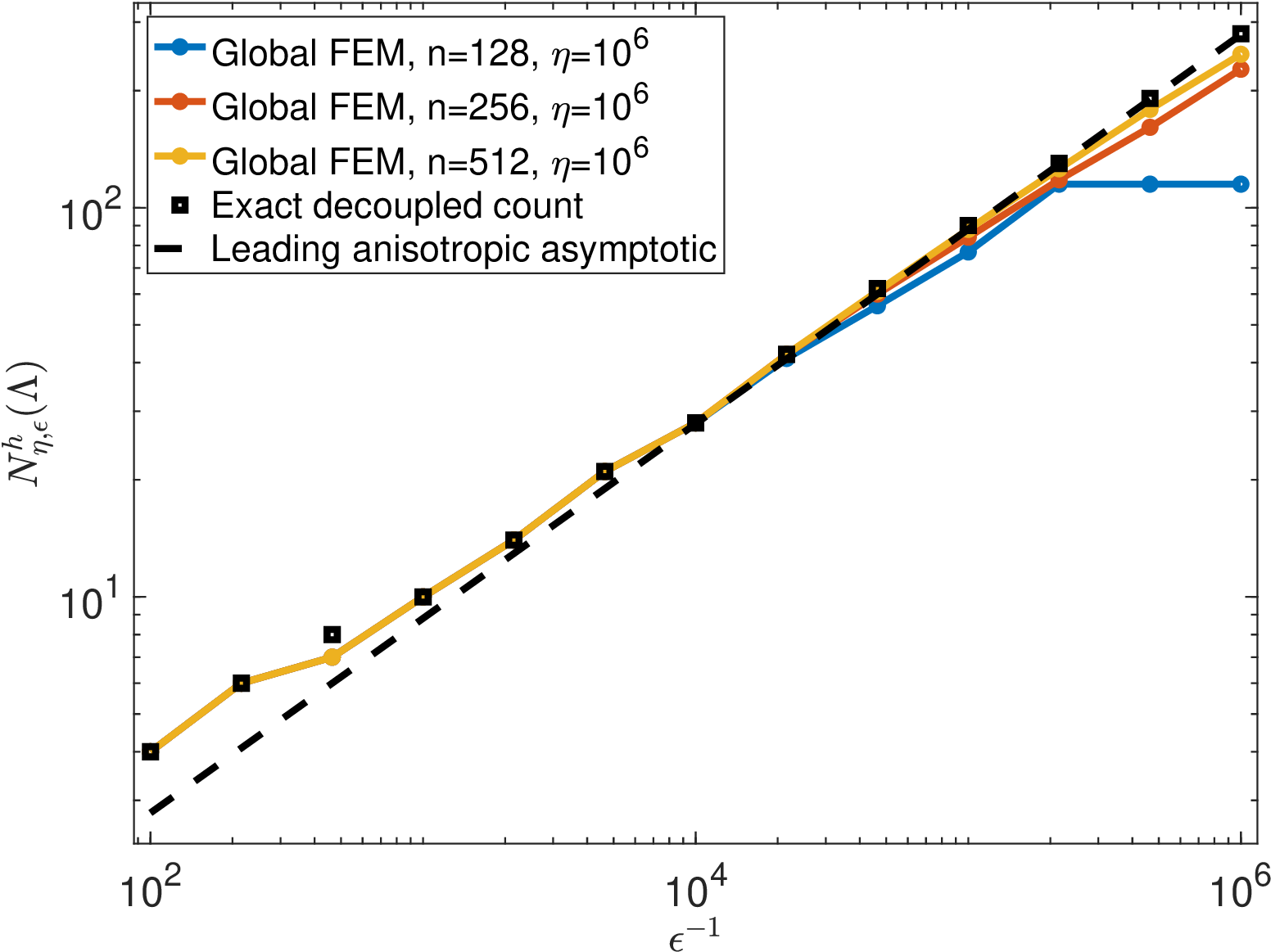}
    \caption{Global finite element eigenvalue counts in $[0,1]$ for $\eta=10^6$ on meshes with $n=128,256,512$, compared with the exact decoupled count given by \eqref{eq:three-inclusion-exact-count} and the leading asymptotic term $(7/(8\pi))\epsilon^{-1/2}$.}
    \label{fig:weyl-count}
\end{figure}

We fix the counting threshold at $\Lambda=1$. The smallest nonzero transverse eigenvalue among the inclusions is $4\pi^2>1$. Since $\Omega_0\subset\Omega$, the first background Dirichlet eigenvalue is at least $2\pi^2>1$. Consequently, only the transverse modes with $j=0$ contribute in this window, and the exact counting formula reduces to
\begin{equation}
N_{\infty,\epsilon}(1) = \sum_{i=1}^{3} \left[ 1+\left\lfloor\frac{L_i}{\pi\sqrt{\epsilon}}\right\rfloor \right].
\label{eq:three-inclusion-exact-count}
\end{equation}
The three local constant modes are included in this count. Since $L_1+L_2+L_3=7/8$, \Cref{thm:anisotropic-weyl} gives
\begin{equation}
N_{\infty,\epsilon}(1) = \frac{7}{8\pi}\epsilon^{-1/2}+O(1), \qquad \epsilon\to0.
\label{eq:three-inclusion-weyl}
\end{equation}
For each value of $\epsilon$, we compute the eigenvalue count of the global finite element problem in $[0,1]$ at $\eta=10^6$ on uniform meshes with $n=128,256,512$. The exact count of the decoupled operator is given by \eqref{eq:three-inclusion-exact-count}. \Cref{fig:weyl-count} compares these numerical counts with the exact integer counts and the leading term in \eqref{eq:three-inclusion-weyl}.

As $\epsilon$ decreases, the global counts show increasingly pronounced mesh dependence. The coarsest-mesh curve levels off, while refinement increases the number of resolved modes and reduces the discrepancy from the decoupled reference. At $\epsilon=10^{-6}$, the global counts are $115$, $227$, and $248$ for $n=128,256,512$, respectively, compared with the exact decoupled count of $280$.

\section{Conclusions}
\label{sec:conclu}
This paper analyzed low-energy eigenvalues in high-contrast PDEs with strongly anisotropic inclusions. Stable extensions of local kernel fields give $O(\eta^{-1})$ upper bounds uniformly in $\epsilon$. Compatible weak fields yield $\lambda^{\eta,\epsilon}_{N+K}\le C_K(\epsilon+\eta^{-1})$ for each fixed $K$, and a discrete lifting argument gives corresponding finite element bounds under the stated stability assumptions. For separable scalar diffusion, we derive an exact counting formula for the local Neumann spectra and establish norm-resolvent decoupling as $\eta\to\infty$ at fixed $\epsilon>0$. The eigenvalue count of the limiting operator grows like $\epsilon^{-1/2}$ at fixed positive energy thresholds as $\epsilon\to0$. The numerical examples illustrate the low-energy modes, mesh-dependent spectral transitions, and scalar decoupling and counting results.

Future work will focus on multiscale approximation and preconditioning methods whose coarse spaces account for material anisotropy and resolved weak modes. 
For periodic media, another challenging and worthwhile direction is to understand the connection between the eigenvalue distribution of the anisotropic problem and the spectrum of the corresponding homogenized operator.

\bibliographystyle{siamplain}
\bibliography{references}
\end{document}

%% file: ex_shared.tex
\usepackage{lipsum}
\usepackage{amsfonts}
\usepackage{graphicx}
\usepackage{epstopdf}
\usepackage{algorithmic}
\ifpdf
  \DeclareGraphicsExtensions{.eps,.pdf,.png,.jpg}
\else
  \DeclareGraphicsExtensions{.eps}
\fi

\newsiamremark{remark}{Remark}
\newsiamremark{hypothesis}{Hypothesis}
\newsiamthm{assumption}{Assumption}

\crefname{hypothesis}{Hypothesis}{Hypotheses}
\newsiamthm{claim}{Claim}
\newsiamremark{fact}{Fact}
\crefname{fact}{Fact}{Facts}

\headers{Anisotropic Multiscale spectrum}{Eric T. Chung, Lijian Jiang and Yajun Wang}

\title{Spectral Gaps and Low-Energy Modes in Anisotropic High-Contrast PDEs
\thanks{Submitted to the editors DATE.
\funding{Eric T. Chung's work is partially supported by the Hong Kong RGC General Research Fund (Project Nos. 14304525 and 14305624). Lijian Jiang acknowledges the support of the NSFC (Project No. 12271408).}}}

\author{
Yajun Wang \thanks{Department of Mathematics, The Chinese University of Hong Kong, Hong Kong SAR (\email{yajunwang@cuhk.edu.hk}).}
\and Eric T. Chung\thanks{Department of Mathematics, The Chinese University of Hong Kong, Hong Kong SAR (\email{eric.t.chung@cuhk.edu.hk}).}
\and Lijian Jiang \thanks{School of Mathematical Sciences, Tongji University, Shanghai, P.R. China (\email{ljjiang@tongji.edu.cn}).}
}

\usepackage{amsopn}
